\documentclass[10pt]{amsart}

\usepackage{verbatim}
\usepackage{eucal,url,amssymb,stmaryrd,
enumerate,amscd,
}

\usepackage{amsfonts}
\usepackage{amsmath,amsthm,amssymb,amscd,enumerate,eucal,url,stmaryrd}
\usepackage{color}

\numberwithin{equation}{section}

\newtheorem{thrm}{Theorem}[section]

\newtheorem{lemma}[thrm]{Lemma}

\newtheorem{prop}[thrm]{Proposition}

\newtheorem{cor}[thrm]{Corollary}

\newtheorem{dfn}[thrm]{Definition}

\newtheorem{rmrk}[thrm]{Remark}

\newtheorem{exm}[thrm]{Example}

\def\frg{{\frak g}}
\def\frh{{\frak h}}

\def\frs{{\frak s}}

\def\fru{{\frak u}}
\def\gc{\frg_\mathbb{C}}
\def\Re{{\frak R}{\frak e}\,}

\def\nilm{\Gamma\backslash G}

\def\db{{\bar{\partial}}}
\def\zzz{{\!\!\!}}

\begin{document}

\title[Balanced Hermitian geometry on 6-dimensional nilmanifolds] {Balanced Hermitian geometry on 6-dimensional nilmanifolds
}

\author{Luis Ugarte}
\address[L. Ugarte]{Departamento de Matem\'aticas\,-\,I.U.M.A.\\
Universidad de Zaragoza\\
Campus Plaza San Francisco\\
50009 Zaragoza, Spain}
\email{ugarte@unizar.es}

\author{Raquel Villacampa}
\address[R. Villacampa]{Centro Universitario de la Defensa\,-\,I.U.M.A.\\
Academia General Militar\\
Crta. de Huesca~s/n.\\
50090 Zaragoza, Spain}
\email{raquelvg@unizar.es}

\maketitle

\begin{abstract}
The invariant balanced Hermitian geometry of
nilmanifolds of dimension 6 is described.
We prove that the (restricted) holonomy group of the
associated Bismut connection reduces to a proper subgroup of SU(3)
if and only if the complex structure is abelian. As an
application we show that if $J$ is abelian then any invariant balanced $J$-Hermitian structure
provides solutions of the Strominger system.
\end{abstract}


\section{Introduction}

\noindent
Given any Hermitian structure $(J,F)$ on
a $2n$-dimensional manifold $M$, Bismut proved in \cite{Bis} the existence of a unique Hermitian connection
with torsion $T$ given by $g(X,T(Y,Z))=JdF(X,Y,Z)=-dF(JX,JY,JZ)$,
$g$ being the associated metric. This torsion connection
will be denoted here by $\nabla$ and the torsion $T$ will be identified with
the 3-form $JdF$. In relation to the Levi-Civita connection
$\nabla^g$ of the Riemannian metric~$g$, the Bismut connection is
determined by $\nabla = \nabla^g +\frac{1}{2}T$.

Since the connection $\nabla$ is Hermitian, its (restricted) holonomy group ${\rm Hol}(\nabla)$ is contained
in the unitary group U($n$). We are interested here in the case when ${\rm Hol}(\nabla)$ is reduced to SU($n$),
and a Hermitian structure satisfying this condition is said to be {\it Calabi-Yau with torsion}.
In dimension six, these structures are related to the Strominger system in heterotic string theory~\cite{Str} and several
constructions of Calabi-Yau with torsion manifolds can be found in~\cite{GP,Grant,GGP,GIP}.

Suppose that $M$ is a nilmanifold, i.e. a compact quotient of a simply-connected nilpotent Lie group $G$ by a lattice,
endowed with an invariant complex structure $J$, i.e. $J$ stems from a complex structure on the Lie algebra $\frg$ of $G$.
It is proved in~\cite{FG} that if there is a $J$-Hermitian metric on $M$ such that the holonomy of the associated Bismut connection is
contained in SU($n$) then there
is an {\it invariant} $J$-Hermitian metric on $M$ which is balanced in the sense of \cite{Mi}. The latter condition means that the associated Lee
1-form $\theta$ vanishes identically or, equivalently, the form $F^{n-1}$ given by the wedge
product of the K\"ahler form $(n-1)$-times is closed.
Moreover, Fino, Parton and Salamon~\cite{FPS} proved that for an invariant Hermitian
structure $(J,F)$ on $M$, the balanced condition is equivalent to the Calabi-Yau with torsion condition.

Our goal in this paper is the study of the invariant balanced Hermitian
geometry of 6-dimensional nilmanifolds, the behaviour
of the holonomy of the associated Bismut connection $\nabla$
and the application to finding solutions of the Strominger system with respect to $\nabla$ in the anomaly
cancellation condition.

In greater detail, the paper is structured as follows.
Section~\ref{complex-sect} is devoted to a detailed description of the invariant balanced Hermitian geometry of 6-dimensional nilmanifolds $M$.
It is proved in~\cite{U} that the Lie algebra $\frg$ underlying $M$ must be isomorphic to $\frh_2$,
$\frh_3$, $\frh_4$, $\frh_5$, $\frh_6$ or
$\frh_{19}^-$. The latter is the only one
for which the complex structure is of non-nilpotent type and its balanced Hermitian geometry is studied in~\cite{UV}.
On the other hand, only $\frh_5$ admits a complex-parallelizable structure $J_0$ and the pair $(\frh_5,J_0)$ corresponds to the well-known Iwasawa manifold.
We put special attention to the balanced Hermitian geometry associated to abelian
complex structures. In the list above, the Lie algebras having abelian complex structures $J$ are $\frh_2$,
$\frh_3$, $\frh_4$ and $\frh_5$, but it turns out
that the $J$-Hermitian metrics on $\frh_2$ or $\frh_4$ are never balanced (see Proposition~\ref{caracterizacion_h3_h5}).
In contrast, any abelian complex structure on $\frh_5$ admits balanced Hermitian
metrics by Corollary~\ref{balanced-h5}. The Lie algebra $\frh_3$
is special, since there exist, up to isomorphism,
two complex structures
but only one of them admits compatible balanced metrics. This Lie algebra corresponds to the product Lie group $H\times \mathbb{R}$, where
$H$ is the 5-dimensional generalized Heisenberg group.

The main result in Section~\ref{complex-sect} is Theorem~\ref{general-balanced-real-basis}, which gives a description of the invariant balanced geometry on 6-dimensional nilmanifolds in terms of a global basis of 1-forms $\{e^1,\ldots,e^6\}$ adapted to the structure $(J,F)$, in the sense that the complex structure $J$ and the fundamental form $F$ express in the canonical way $Je^1=-e^2,\ Je^3=-e^4,\ Je^5=-e^6$ and $F=e^{12}+e^{34}+e^{56}$.

In Section~\ref{deform} we study on nilmanifolds the weak $\partial\db$-lemma recently introduced by Fu and Yau \cite{FY2}
in relation to deformations of balanced metrics. It is proved in \cite{FY2} that given a compact complex $n$-dimensional manifold $M$ with a balanced metric, if along a small deformation $M_{\lambda}$ of $M$ the $(n-1,n)$-th weak $\partial\db$-lemma is satisfied then
there exists a balanced metric on $M_{\lambda}$ for sufficiently small $\lambda$.
If $M$ is a nilmanifold endowed with an invariant complex structure~$J$, using the symmetrization process
and the results on the Dolbeault cohomology of $(M,J)$ obtained by Rollenske in \cite{R}, we show when the
weak $\partial\db$-lemma on $(M,J)$ is reduced to the study of the weak $\partial\db$-lemma at the Lie algebra level.
In particular, if the complex structure $J$ is abelian then the $(n-1,n)$-th weak $\partial\db$-lemma is always satisfied.
The general behaviour in dimension 6 with respect to the weak $\partial\db$-lemma in the presence of balanced
structures is also described. As an application, we give an explicit deformation $I_{\lambda}$ of an abelian complex structure $I_0$ on a nilmanifold
associated to $\frh_5$ having compatible balanced metric such that
the $(2,3)$-th weak $\partial\db$-lemma only holds for $\lambda=0$ but $I_{\lambda}$ admits balanced metric for any $\lambda$,
which shows that the weak $\partial\db$-lemma is not a necessary condition for the existence of balanced metric along deformation
of the complex structure.

In Section~\ref{Bismut-hol}, using the description given in Theorem~\ref{general-balanced-real-basis}, we determine the
(restricted) holonomy group of the Bismut connection $\nabla$ for any invariant balanced Hermitian structure $(J,F)$.
We prove in Theorem~\ref{Main1} that ${\rm Hol}(\nabla)={\rm SU(3)}$ if and only if $J$ is not abelian.
In the abelian case it is rather straightforward to show that the holonomy reduces to a subgroup of SU(2) (see Remark~\ref{paralelismo} for details),
so the main effort in proving the result is to verify that ${\rm Hol}(\nabla)$ actually equals SU(3) in all the remaining cases,
for which we study the behaviour of the curvature endomorphisms
of $\nabla$ and their covariant derivatives of any order.
Since the Bismut connection $\nabla$ depends on the pair $(J,F)$, it is a surprising fact that for 6-dimensional nilmanifolds the behaviour of
${\rm Hol}(\nabla)$ is determined by the complex structure $J$
and it does not depend on the balanced compatible metric.
It is also proved that if the complex structure is abelian then ${\rm Hol}(\nabla)={\rm SU(2)}$
if and only if the underlying Lie algebra is $\frh_5$, i.e.
the Lie algebra corresponding to the Iwasawa manifold.
In Example~\ref{abelian-solvable} we present an abelian complex structure $J$
on a compact solvmanifold of dimension 6 admitting an invariant balanced $J$-Hermitian metric such that the holonomy
group of its associated Bismut connection equals SU(3), so Theorem~\ref{Main1} cannot be extended to solvmanifolds.

As an application, in the last section we look for solutions of the Strominger system~\cite{Str} in the class of
invariant balanced Hermitian 6-dimensional nilmanifolds, which forces the dilaton function to be constant (see
equations (a)-(d) in Section~\ref{het-constant-dil} for details).
In the Strominger system, the anomaly cancellation condition can be solved for different choices of metric connection and the
physical validity of the corresponding solutions is studied in~\cite{BBFTY}.
In this context, it is relevant a recent result by Ivanov~\cite{I} asserting that a solution of the Strominger system provides also
a solution of the heterotic equations of motion if and only if the metric connection is an instanton.
In~\cite{FY} Fu and Yau consider the Chern connection $\nabla^c$ and prove the existence of solutions with non-constant dilaton
on a Hermitian non-K\"ahler manifold given as a $\mathbb{T}^2$-bundle over a $K3$ surface.
In the case of constant dilaton, explicit solutions are given in \cite{FIUV} based on $\frh_2,\ldots,\frh_6$ and
$\frh_{19}^-$ for different choices of connection, including
the Bismut and Chern connections, and in addition solutions of the heterotic equations of motion were found for $\frh_3$.
A recent solution based on $\frh_3$ with respect to another metric connection is obtained by Grantcharov~\cite{Grant}.
In \cite{UV} it is proved that for any invariant complex structure $J$ on a nilmanifold $N$ with $\frh_{19}^-$ as underlying Lie algebra,
the compact complex manifold $(N,J)$
admits solutions, with constant dilaton and non-flat instanton, of the Strominger system satisfying the anomaly
cancellation condition with respect to the Chern connection $\nabla^c$.

In Theorems~\ref{solutions-h3} and~\ref{solutions-h5} we prove that any abelian complex structure provides solutions of the Strominger system. More concretely, let $M$ be a nilmanifold endowed with an invariant balanced
Hermitian structure $(J,F)$. If $J$ is abelian, then
there is an invariant non-flat SU(3)-instanton solving the Strominger system with respect to the Bismut connection in the anomaly
cancellation condition.
Moreover, any such solution solves in addition the heterotic equations of motion if and only if
$\frh_3$ is the Lie algebra underlying $M$.
Finally, in Section~\ref{more-solutions} more solutions for non-abelian complex structures are given.
When $J$ is of non-nilpotent type, we prove the existence of a non-flat instanton solving at the same time the Strominger systems for
the Bismut and the Chern connection (see Proposition~\ref{simultaneous}).

\section{Invariant complex structures on 6-dimensional nilmanifolds 
and compatible balanced metrics}\label{complex-sect}

\noindent Let $M$ be a nilmanifold of even dimension, i.e. a compact quotient of a
simply-connected nilpotent Lie group $G$ by a lattice $\Gamma$ of
maximal rank. Any left-invariant complex structure on $G$ descends
to $M$ in a natural way, so a source (possibly empty) of complex
structures on $M$ is given by the endomorphisms $J\colon\frg
\longrightarrow \frg$ of the Lie algebra $\frg$ of $G$ such that
$J^2=-{\rm Id}$ satisfying the ``Nijenhuis condition''
$$
[JX,JY]=J[JX,Y]+J[X,JY]+[X,Y],
$$
for any $X,Y\in \frg$. We shall refer to any such an endomorphism as
a {\it complex structure} on the Lie algebra~$\frg$.

Associated to a complex structure $J$, there exists an ascending
series $\{\frg_l^J\}_{l\geq 0}$ of the Lie algebra defined
inductively by
$$
\frg_0^J = \{0\} \ , \qquad  \frg_l^J = \{ X \in \frg \, \mid \, [X,
\frg] \subseteq \frg_{l-1}^J \ \ \mbox{\rm and} \ \ [JX, \frg]
\subseteq \frg_{l-1}^J \}\ , \quad l\geq 1.
$$
For any $l\geq 0$, the term $\frg_l^J$ is a $J$-invariant ideal of
$\frg$ which is contained in the term $\frg_l=\{X\in\frg\mid
[X,\frg]\subseteq \frg_{l-1}\}$ of the usual ascending central
series of $\frg$. But whereas $\{\frg_l\}_{l\geq 0}$ always reaches the whole Lie
algebra when $\frg$ is nilpotent, the series
$\{\frg_l^J\}_{l\geq 0}$ can stabilize in a proper $J$-ideal of
$\frg$. This motivates the following terminology: if $\frg_l^J=\frg$
for some $l$ then the complex structure $J$ is called {\it
nilpotent}~\cite{CFGU2}; otherwise, we shall say that $J$ is {\it
non-nilpotent}.

Well-known particular classes of nilpotent complex structures are the
{\it complex-parallelizable} structures, for which $[JX,Y]=J[X,Y]$, and the {\it abelian} structures, which satisfy the condition $[JX,JY]=[X,Y]$.
A Lie algebra $\frg$ has a complex-parallelizable structure if and only if $\frg$ can be endowed with a complex Lie algebra structure.

\begin{dfn}\label{def-non-nilp}
{\rm A {\it nilpotent} (resp. {\it non-nilpotent}) complex structure
on a nilmanifold $M$ is a complex structure on $M$ coming from a
nilpotent (resp. non-nilpotent) complex structure $J$ on the
underlying Lie algebra $\frg$.}
\end{dfn}

Let us denote by $\gc$ the complexification of $\frg$ and by $\gc^*$
its dual. Given an endomorphism $J\colon \frg \longrightarrow \frg$
such that $J^2=-{\rm Id}$, we denote by $\frg^{1,0}$ and
$\frg^{0,1}$ the eigenspaces corresponding to the eigenvalues $\pm i$ of $J$ as an
endomorphism of~$\gc^*$, respectively. The decomposition
$\gc^*=\frg^{1,0}\oplus\frg^{0,1}$ induces a natural bigraduation on
the complexified exterior algebra $\bigwedge^* \,\gc^* =\oplus_{p,q}
\bigwedge^{p,q}(\frg^*)=\oplus_{p,q} \bigwedge^p(\frg^{1,0})\otimes
\bigwedge^q(\frg^{0,1})$. If $d$ denotes the usual
Chevalley-Eilenberg differential of the Lie algebra, we shall also
denote by $d$ its extension to the complexified exterior algebra,
i.e. $d\colon \bigwedge^* \gc^* \longrightarrow \bigwedge^{*+1}
\gc^*$. It is well-known that the endomorphism $J$ is a complex
structure if and only if $d(\frg^{1,0})\subset
\bigwedge^{2,0}(\frg^*)\oplus \bigwedge^{1,1}(\frg^*)$. In the case
of nilpotent Lie algebras $\frg$, Salamon proves in~\cite{S} the
following equivalent condition for the endomorphism $J$ to be a
complex structure: $J$ is a complex structure on $\frg$ if and only
if $\frg^{1,0}$ has a basis $\{\omega^j\}_{j=1}^n$ such that
$d\omega^1=0$ and
$$
d \omega^{j} \in \mathcal{I} (\omega^1,\ldots,\omega^{j-1}), \quad
\mbox{ for } j=2,\ldots,n ,
$$
where $\mathcal{I} (\omega^1,\ldots,\omega^{j-1})$ is the ideal in
$\bigwedge\phantom{\!}^* \,\gc^*$ generated by
$\{\omega^1,\ldots,\omega^{j-1}\}$. From now on, we shall denote
$\omega^{j}\wedge\omega^k$ and $\omega^{j}\wedge\overline{\omega^k}$
simply by $\omega^{jk}$ and $\omega^{j\bar{k}}$, respectively.

A complex structure $J$ is nilpotent if and only if there is a basis $\{\omega^j\}_{j=1}^n$
for~$\frg^{1,0}$ satisfying $d\omega^1=0$ and
$$
d \omega^j \in \bigwedge\phantom{\!\!}^2 \,\langle
\omega^1,\ldots,\omega^{j-1},
\omega^{\overline{1}},\ldots,\omega^{\overline{j-1}} \rangle, \quad
\mbox{ for } j=2,\ldots,n .
$$
Abelian complex structures satisfy in addition that $d(\frg^{1,0}) \subset \bigwedge^{1,1}(\frg^*)$, and
they are characterized by the fact that the complex Lie algebra $\frg^{1,0}$ is abelian.
Finally, a nilpotent complex structure is complex-parallelizable if and only if $d(\frg^{1,0}) \subset \bigwedge^{2,0}(\frg^*)$.

\smallskip

Now, let $\frg$ be a Lie algebra of dimension 6. A {\it Hermitian
structure} on $\frg$ is a pair $(J,g)$, where $J$ is a complex
structure on $\frg$ and $g$ is an inner product on $\frg$ compatible
with $J$ in the usual sense, i.e. $g(\cdot,\cdot)=g(J\cdot,J\cdot)$.
The associated {\it fundamental form} $F\in \bigwedge^2 \frg^*$ is
defined by $F(X,Y)=g(X,JY)$ and expresses in terms of any basis
$\{\omega^j\}_{j=1}^3$, of type (1,0) with respect to $J$, by
\begin{equation}\label{2forma}
2\, F=i (r^2\omega^{1\bar{1}} + s^2\omega^{2\bar{2}} +
t^2\omega^{3\bar{3}})
+u\,\omega^{1\bar{2}}-\bar{u}\,\omega^{2\bar{1}} +
v\,\omega^{2\bar{3}}-\bar{v}\,\omega^{3\bar{2}}+
z\,\omega^{1\bar{3}}-\bar{z}\,\omega^{3\bar{1}},
\end{equation}
for some $r,s,t\in\mathbb{R}$ and $u,v,z\in\mathbb{C}$. Since we are
using the convention $(J\alpha)(X)=-\alpha(JX)$ for $X\in\frg$ and
$\alpha\in\frg^*$, the inner product $g$ is given by
\begin{equation}\label{metric}
g= r^2\,\omega^1\omega^{\bar{1}} + s^2\, \omega^2\omega^{\bar{2}} +
t^2\, \omega^3\omega^{\bar{3}} - \frac{i}{2}(u\,
\omega^1\omega^{\bar{2}} - \bar{u}\, \omega^2\omega^{\bar{1}} +v\,
\omega^2\omega^{\bar{3}} - \bar{v}\, \omega^3\omega^{\bar{2}}+ z\,
\omega^1\omega^{\bar{3}} -\bar{z}\, \omega^3\omega^{\bar{1}}).
\end{equation}
Here $\omega^j\omega^{\bar{k}}=\frac12
(\omega^j\otimes\omega^{\bar{k}} + \omega^{\bar{k}}\otimes\omega^j)$
denotes the symmetric product of $\omega^j$ and $\omega^{\bar{k}}$.
Notice that the positive definiteness of $g$ implies that the
coefficients $r^2,\,s^2,\,t^2$ are non-zero real numbers and
$u,\,v,\,z\in \mathbb{C}$ satisfy $r^2s^2>|u|^2,$
$s^2t^2>|v|^2,\,r^2t^2>|z|^2$ and $r^2s^2t^2 + 2\,\Re(i\bar u\bar
vz)>t^2|u|^2 + r^2|v|^2 + s^2|z|^2$.

Fixed $J$, since $g$ and $F$ are mutually determined by each other,
we shall also denote the Hermitian structure $(J,g)$ by the pair
$(J,F)$. Recall that the Hermitian structure $(J,F)$ is said to be
{\it balanced} if $F^2$ is a closed form or, equivalently, $F\wedge
dF=0$.

\begin{dfn}\label{def-inv-balanced}
{\rm An {\it invariant balanced Hermitian structure}
on a nilmanifold $M$ is a balanced Hermitian structure on $M$ coming from a balanced
Hermitian structure $(J,F)$ on the Lie algebra $\frg$ underlying $M$.}
\end{dfn}

The goal of this section is to explicitly describe the invariant balanced Hermitian geometry of 6-dimensional nilmanifolds.

\begin{prop}\label{general-balanced}
Let $(J,F)$ be a balanced Hermitian structure on a $6$-dimensional
(non-abelian) nilpotent Lie algebra $\frg$.
\begin{enumerate}
\item[{\rm (i)}] If $J$ is a complex-parallelizable structure, then
$\frg^*_{\mathbb{C}}$ has a basis
$\{\omega^j,\omega^{\bar{j}}\}_{j=1}^3$ such that
\begin{equation}\label{general_complex_parallelizable}
d\omega^1=d\omega^2=0,\quad
d\omega^3=\omega^{12},\end{equation}
and any $J$-Hermitian structure is balanced.
\item[{\rm (ii)}]
If the complex structure $J$ is nilpotent (but not complex
parallelizable) then there exists a basis
$\{\omega^j,\omega^{\bar{j}}\}_{j=1}^3$ for $\frg^*_{\mathbb{C}}$
satisfying
\begin{equation}\label{J-nilp}
d\omega^1 = 0,\quad\quad d\omega^2 = 0,\quad\quad  d\omega^3 =
\rho\, \omega^{12} + \omega^{1\bar{1}} + b^2\, \omega^{1\bar{2}} +
(x+y i)\, \omega^{2\bar{2}},
\end{equation} where $\rho=0,1$ and $b,x,y\in \mathbb{R}$, such that the fundamental form $F$ expresses as
\begin{equation}\label{F-nilp}
2\,F=i\,(\omega^{1\bar1}+s^2\,\omega^{2\bar2}+t^2\,\omega^{3\bar3})+u\,\omega^{1\bar2}-\bar
u\,\omega^{2\bar1},
\end{equation}
where $s^2>|u|^2$ and $t^2>0$ satisfy
\begin{equation}\label{balanced-condition-nilp}
s^2+x+y\, i = \bar u\, b^2 i.
\end{equation}
\item[{\rm (iii)}]
If $J$ is non-nilpotent then there is a basis
$\{\omega^j,\omega^{\bar{j}}\}_{j=1}^3$ for $\frg^*_{\mathbb{C}}$
satisfying
\begin{equation}\label{J-non-nilp}
d\omega^1 = 0,\quad\quad d\omega^2 = \omega^{13} +
\omega^{1\bar{3}}, \quad\quad d\omega^3 = \pm i\, (\omega^{1\bar{2}}
- \omega^{2\bar{1}}) ,
\end{equation}
and the fundamental form $F$ expresses as
\begin{equation}\label{F-non-nilp}
2\,F=i(r^2\,\omega^{1\bar1} + s^2\,\omega^{2\bar2} +
t^2\,\omega^{3\bar3}) + v\,\omega^{2\bar3} -\bar v\,\omega^{3\bar2},
\end{equation}
where $r^2>0$ and $s^2t^2>|v|^2$.
\end{enumerate}
\end{prop}

\begin{proof}
The assertion (i) is well-known and the proof of (iii) is given in \cite{UV}, so it remains to prove (ii).
By \cite{U}, if $\frg$ has a balanced structure compatible with a nilpotent complex structure $J$ then the Lie algebra is 2-step nilpotent.
Moreover \cite{U}, for any nilpotent (not complex-parallelizable) complex structure $J$
on a 2-step nilpotent Lie algebra $\frg$, there exists a
$(1,0)$-basis $\{\omega'^j\}_{j=1}^3$ such that
$$
d\omega'^1 = 0,\quad d\omega'^2 = 0,\quad d\omega'^3 = \rho\,
\omega'^{12} + \omega'^{1\bar{1}} + B\, \omega'^{1\bar{2}} + D\,
\omega'^{2\bar{2}},
$$
where $\rho=0,1$ and $B,\,D\in \mathbb{C}$. Suppose $B\neq 0$ and let
$\zeta$ be any non-zero solution of the equation $\bar
\zeta\frac{B}{|B|}=\zeta$. We can choose $\zeta$ of modulo 1 and
with respect to the basis
$\{\omega^1=\zeta\,\omega'^1,\,\omega^2=\bar
\zeta\,\omega'^2,\,\omega^3=\omega'^3\}$ one gets (\ref{J-nilp})
with coefficient $b^2=|B|>0$.

Now, consider (\ref{J-nilp}) expressed in terms of a basis
$\{\tau^1,\,\tau^2,\,\tau^3\}$ and a general structure
(\ref{2forma}). The (1,0)-basis given by
$\{\sigma^1=\tau^1,\,\sigma^2=\tau^2,\,\sigma^3=-\frac{iz}{t^2}\,\tau^1-\frac{iv}{t^2}\,\tau^2+\tau^3\}$
preserves the complex equations (\ref{J-nilp}) and the fundamental
form $F$ expresses in terms of $\{\sigma^1,\sigma^2,\sigma^3\}$ as
$$2\, F=i\,(r'^2\,\sigma^{1\bar 1} + s'^2\,\sigma^{2\bar2} +
t'^2\,\sigma^{3\bar3}) + u'\,\sigma^{1\bar2} - \bar
u'\,\sigma^{2\bar1},$$ with new metric coefficients
$r'^2=r^2-\frac{|z|^2}{t^2},\, s'^2=s^2-\frac{|v|^2}{t^2},\,
t'^2=t^2$ and $u'=u-\frac{i\bar vz}{t^2}$. Moreover, we get
(\ref{F-nilp}) after normalizing the coefficient $r'^2$ by
considering
$\{\omega^1=r'\,\sigma^1,\,\omega^2=r'\,\sigma^2,\,\omega^3=r'^2\,\sigma^3\}$,
which also preserves the equations (\ref{J-nilp}).

A direct calculation shows that a Hermitian structure given by
(\ref{J-nilp}) and (\ref{F-nilp}) satisfies
\begin{eqnarray*}
4\,F\wedge dF&=&t^2\,(s^2+x+y i- \bar u\, b^2
i)\,\omega^{12\bar1\bar2\bar3} + t^2\,(s^2+x-y i + u\, b^2
i)\,\omega^{123\bar1\bar2}.
\end{eqnarray*}
Therefore, the Hermitian structure $(J,F)$ is balanced if and only
if $s^2+x+y i= \bar u\, b^2 i$. This completes the proof of (ii).
\end{proof}

The nilpotent Lie algebras $\frg$ admitting balanced Hermitian structure are
classified in \cite{U}. They are: $\frh_2=(0,0,0,0,12,34)$,
$\frh_3=(0,0,0,0,0,12+34)$, $\frh_4=(0,0,0,0,12,14+23)$,
$\frh_5=(0,0,0,0,13+42,14+23)$, $\frh_6=(0,0,0,0,12,13)$ and
$\frh_{19}^-=(0,0,0,12,23,14-35)$. The latter is the only one
corresponding to the case (iii) in the previous proposition,
i.e. the complex structure is non-nilpotent; moreover, up to isomorphism there exist only
two complex structures $J_0^{\pm}$ on the Lie algebra $\frh_{19}^-$ \cite{UV}, which correspond
to the $\pm$-sign in \eqref{J-non-nilp}, respectively.
On the other hand, it is well-known that only $\frh_5$ admits a complex-parallelizable structure, and
we shall denote by $J_0$ the structure given by~\eqref{general_complex_parallelizable}.
Notice that the pair $(\frh_5,J_0)$ corresponds to the Iwasawa manifold.

We shall use the following result, which gives a classification of
the Lie algebras underlying the structure equations~(\ref{J-nilp})
depending on the values of the quadruplet $(\rho,b^2,x,y)$.

\begin{lemma}\label{underlyingLiealgebras} \cite[Proposition 13]{U}
Let $J$ be a complex structure on a nilpotent Lie algebra $\frg$ given by
(\ref{J-nilp}). Then:
\begin{itemize}
\item[(i)] If $b^2=\rho$, then the Lie algebra $\frg$ is isomorphic
to:
\begin{itemize}
\item[(i.1)]  $\frh_2$ for $y\neq 0$;
\item[(i.2)] $\frh_3$ for $\rho=y=0$ and $x\neq 0$;
\item[(i.3)] $\frh_4$ for $\rho=1,\, y=0$ and $x\neq 0$;
\item[(i.4)] $\frh_6$ for $\rho=1$ and $x=y=0$.
\end{itemize}
\item[(ii)] If $b^2\neq \rho$, then the Lie algebra $\frg$ is
isomorphic to:
\begin{itemize}
\item[(ii.1)] $\frh_2$ for $4y^2>(\rho-b^4)(4x+\rho-b^4)$;
\item[(ii.2)] $\frh_4$ for $4y^2=(\rho-b^4)(4x+\rho-b^4)$;
\item[(ii.3)] $\frh_5$ for $4y^2<(\rho-b^4)(4x+\rho-b^4)$.
\end{itemize}
\end{itemize}
\end{lemma}

We have omitted the case $b^2=\rho=x=y=0$, which corresponds to the
Lie algebra $\frh_8=(0,0,0,0,0,12)$, because it does not admit any balanced
Hermitian structure.

From now on we shall concentrate mainly in the case when the complex
structure is nilpotent, because the non-nilpotent case is studied in detail
in \cite{UV}, although we will return to it in Proposition~\ref{reduc-hol-caso-non-nilp}.
In view of Proposition \ref{general-balanced}~(ii), if we denote by $u_1$ and $u_2$ the
real and imaginary parts of $u$, i.e. $u=u_1+u_2i$,
then the balanced condition (\ref{balanced-condition-nilp}) reads as
$x=u_2b^2-s^2$ and $y=u_1b^2$. Therefore, the balanced Hermitian
structures $(J,F)$, $J$ nilpotent, are parametrized by $\rho=0,1$
and a 5-tuple $(b,u_1,u_2,s^2,t^2)\in \mathbb{R}^3\times
\mathbb{R}^+\times \mathbb{R}^+$ satisfying $s^2-u_1^2-u_2^2>0$, in
the sense that the complex structure $J$ is given by
$$d\omega^1 = d\omega^2 = 0,\quad\quad  d\omega^3 = \rho\,
\omega^{12} + \omega^{1\bar{1}} + b^2\, \omega^{1\bar{2}} +
(u_2b^2-s^2+u_1b^2 i)\, \omega^{2\bar{2}},$$ and the fundamental
form $F$ expresses as
$$
2\,F=i(\omega^{1\bar1}+ s^2\, \omega^{2\bar2}+
t^2\,\omega^{3\bar3}) + (u_1+u_2i)\omega^{1\bar2}-(u_1-u_2i)\omega^{2\bar1}.
$$

We recall that a Hermitian structure $(J,F)$ on $\frg$ is said to be
{\it equivalent} to a Hermitian structure $(J',F')$ on $\frg'$ if
there is an isomorphism $A\colon \frg \longrightarrow \frg'$ of Lie
algebras such that $J'A=AJ$ and $F=A^*F'$.

\begin{lemma}\label{h5_complex_parallelizable}
On $(\frh_5,J_0)$ any (balanced) Hermitian structure is equivalent to
one and only one structure in the 1-parameter family
$2F=i\,(\omega^{1\bar1}+\omega^{2\bar2}+t^2\,\omega^{3\bar3}),\,
t\neq0.$
\end{lemma}

\begin{proof}
Let us consider a generic $J_0$-Hermitian structure
$$
2F=i(r^2\omega^{1\bar1}+s^2\omega^{2\bar2}+t^2\omega^{3\bar3})+ u\,\omega^{1\bar2} - \bar u\,\omega^{2\bar1}
+ v\,\omega^{2\bar3} - \bar v\,\omega^{3\bar2} + z\,\omega^{1\bar3}
- \bar z\,\omega^{3\bar1},
$$
where $\{\omega^j, \omega^{\bar j}\}_{j=1}^3$ is the basis satisfying \eqref{general_complex_parallelizable}.
We define the new $(1,0)$-basis
$\{\mu^1=\omega^1,\,\mu^2=\omega^2,\,\mu^3=-\frac{iz}{t^2}\,\omega^1-\frac{iv}{t^2}\,\omega^2+\omega^3\}$
in order to get $F$ expressed as
$$F=i\,(r'^2\mu^{1\bar1}+s'^2\mu^{2\bar2}+t'^2\mu^{3\bar3})+u'\,\mu^{1\bar2}-\overline{u'}\,\mu^{2\bar1}.
$$
Now, the basis
$\{\sigma^1=i\sqrt{\frac{r'^2s'^2-|u'|^2}{s'^2}}\,\mu^1,\,\sigma^2=-\frac{u'i}{s'}\,\mu^1
+ s'\,\mu^2,\,\sigma^3=i\sqrt{r'^2s'^2-|u'|^2}\,\mu^3\}$ satisfies \eqref{general_complex_parallelizable} and the
$J_0$-Hermitian form expresses as $2F_{t''}=i\,(\sigma^{1\bar1}+\sigma^{2\bar2}+t''^2\,\sigma^{3\bar3})$.
Finally, it is
straightforward to verify that two such Hermitian structures $(J_0,F_{t''_1})$ and $(J_0,F_{t''_2})$ are
equivalent if and only if $(t''_1)^2=(t''_2)^2$.
\end{proof}

For general nilpotent structures the situation is rather complicated, however we have the
following partial classification result:

\begin{lemma}\label{equivalence-u=0}
Let $J$ be a nilpotent complex structure given by \eqref{J-nilp} and let
$F_{t}$ and $F_{t'}$ be two balanced $J$-Hermitian structures given by~\eqref{F-nilp}
and~\eqref{balanced-condition-nilp} with $u=0$. Then, $(J,F_{t})$ and $(J,F_{t'})$ are equivalent if and only if
$t^2=t'^2$.
\end{lemma}

\begin{proof}
First of all observe that if we fix the complex structure and
$u=0$, then $y=0$, $s^2=-x$ and therefore the only free parameter is the metric coefficient $t^2$.
Let us consider two $(1,0)$-bases $\{\omega^j\}_{j=1}^3$ and
$\{\sigma^j\}_{j=1}^3$ satisfying~\eqref{J-nilp} for the same complex parameters $\rho, b^2, x$ ($y=0$) and
let $F_{t}$ and $F_{t'}$ be two balanced $J$-Hermitian structures given by
$$F_{t}=\frac i2(\omega^{1\bar1}+s^2\omega^{2\bar2}+t^2\omega^{3\bar3}),\quad
F_{t'}=\frac
i2(\sigma^{1\bar1}+s^2\sigma^{2\bar2}+t'^2\sigma^{3\bar3}),\textrm{
where } s^2=-x.$$
Suppose that there exists an equivalence $A\colon \frg \longrightarrow \frg$ between
the two Hermitian structures $(J,F_{t})$ and $(J,F_{t'})$. Since the linear isomorphism
$A^*\colon \frg^*\longrightarrow \frg^*$ commutes with the Chevalley-Eilenberg
differential and the extension of $A^*$ to the complexified
exterior algebra preserves the bigraduation induced by $J$, then
$$\sigma^j=a_{j1}\,\omega^1+a_{j2}\,\omega^2+a_{j3}\,\omega^3,\quad\quad j=1,2,3,$$
where $(a_{jk})\in$ GL($3,\mathbb{C}$) satisfying
$d\sigma^j=a_{j1}\,d\omega^1+a_{j2}\,d\omega^2+a_{j3}\,d\omega^3$, for $j=1,2,3$. This is
equivalent to the conditions
\begin{equation}\label{grupo1}\begin{array}{rcl}
0&\!\!=&\!\!a_{13}=a_{23},\\
\rho\,a_{33}&\!\!=&\!\!\rho\,(a_{11}a_{22}-a_{12}a_{21}),\\
a_{33}&\!\!=&\!\!|a_{11}|^2 + b^2\,a_{11}\overline{a_{21}} - s^2\,|a_{21}|^2,\\
-s^2\,a_{33}&\!\!=&\!\!|a_{12}|^2 + b^2\,a_{12}\overline{a_{22}} -
s^2\,|a_{22}|^2,\\
b^2\,a_{33}&\!\!=&\!\!a_{11}\overline{a_{12}} +
b^2\,a_{11}\overline{a_{22}} -
s^2\,a_{21}\overline{a_{22}},\\
0&\!\!=&\!\!a_{12}\overline{a_{11}} + b^2\,a_{12}\overline{a_{21}} -
s^2\,a_{22}\overline{a_{21}}.
\end{array}\end{equation}
Moreover, the condition $F_t=A^*F_{t'}$ implies that the coefficients
$a_{jk}$ must satisfy the following extra equations
\begin{equation}\label{grupo2}\begin{array}{rcl}
0&\!\!=&\!\!a_{31}=a_{32},\\
1&\!\!=&\!\!|a_{11}|^2+s^2\,|a_{21}|^2,\\
s^2&\!\!=&\!\!|a_{12}|^2+s^2\,|a_{22}|^2,\\
0&\!\!=&\!\!a_{11}\overline{a_{12}} + s^2 a_{21}\overline{a_{22}},\\
t^2&\!\!=&\!\!t'^2\,|a_{33}|^2.
\end{array}\end{equation}
Combining the last equation in~\eqref{grupo1} with the fourth equation in~\eqref{grupo2}
one gets $a_{12}(2\,\overline{a_{11}} +
b^2\,\overline{a_{21}})=0$. We have two possibilities:
if $a_{12}=0$, then it follows from~\eqref{grupo2} that $a_{21}=0$
and $|a_{22}|=1$, and thus the fourth equation in~\eqref{grupo1}
implies that $a_{33}=1$ and therefore $t'^2=t^2$;
on the other hand, if $2\,\overline{a_{11}} +
b^2\,\overline{a_{21}}=0$, then the second equation in~\eqref{grupo2}
reads as $1=(\frac{b^4}{4} + s^2)\,|a_{21}|^2$ and the third
equation in~\eqref{grupo1} expresses as $a_{33}=-(\frac{b^4}{4} +
s^2)\,|a_{21}|^2$, which implies $|a_{33}|=1$ and again $t'^2=t^2$.
\end{proof}

\subsection{Abelian complex structures}

Abelian complex structures correspond to the coefficient $\rho=0$ in the
general structure equations~\eqref{J-nilp} and they are
characterized by the condition $[JX,JY]=[X,Y]$, which is equivalent to say that
$d(\frg^{1,0}) \subset \bigwedge^{1,1}(\frg^*)$. Next we study the balanced
Hermitian geometry associated to abelian
complex structures. First we observe that in the case of
abelian structures, Proposition~\ref{general-balanced}~(ii) can be improved
in the sense that coefficient $b^2$ can be reduced to take the value 0
or~1.

\begin{lemma}\label{propreducednilpotentAbelian}
Let $J$ be an abelian complex structure on a $6$-dimensional
$2$-step nilpotent Lie algebra $\frg$.  Then, there is a $(1,0)$-basis
$\{\omega^j\}_{j=1}^3$ satisfying
\begin{equation}\label{reducednilpotentAbelian}
d\omega^1 = d\omega^2 = 0,\quad\quad
d\omega^3 = \omega^{1\bar{1}} +\delta\,
\omega^{1\bar{2}} + (x+y\, i)\, \omega^{2\bar{2}},
\end{equation}
where $\delta=0,1$ and $x,y\in \mathbb{R}$.
\end{lemma}

\begin{proof}
By Proposition~\ref{general-balanced} there is a ($1,0$)-basis
$\{\omega'^j\}^3_{j=1}$ satisfying (\ref{J-nilp}) with
$\rho=0$. Suppose $b^2\neq 0$. We can normalize the coefficient $b^2$ to
be 1 by considering the new basis
$\{\omega^1=\omega'^1,\,\omega^2=b^2\,\omega'^2,\,\omega^3=\omega'^3\}$.
\end{proof}

In the abelian complex case, given a Hermitian structure~\eqref{F-nilp}, the balanced
condition~\eqref{balanced-condition-nilp} reads as
\begin{equation}\label{balanced-condition-abelian}
s^2+x = \delta\, u_2, \quad\quad y=\delta\, u_1.
\end{equation}
Now, we can combine Lemmas~\ref{underlyingLiealgebras} and~\ref{propreducednilpotentAbelian}
together with~\eqref{balanced-condition-abelian} to derive the following classification of nilpotent Lie algebras
admitting abelian complex structures with
balanced compatible metric.

\begin{prop}\label{caracterizacion_h3_h5}
Let $\frg$ be a $6$-dimensional $2$-step nilpotent Lie algebra endowed with an abelian
complex structure $J$ given by~\eqref{reducednilpotentAbelian}. Suppose that
$J$ admits a balanced $J$-Hermitian metric. Then $\frg$ is
isomorphic to $\frh_3$ when $\delta=0$ or $\frh_5$ when $\delta=1$.
\end{prop}

\begin{proof}
The proof follows from Lemma~\ref{underlyingLiealgebras} for the case $\rho=0$. If $\delta=0$ then $y=0$
by~\eqref{balanced-condition-abelian} and the Lie algebra is isomorphic to $\frh_3$. In the case
$\delta=1$, the possibilities for $\frg$ are: $$\frh_2,\, \textrm{ if }
4y^2>1-4x;\quad \frh_4,\, \textrm{ if } 4y^2=1-4x;\quad \frh_5,\,
\textrm{ if } 4y^2<1-4x.$$
Since $s^2>|u|^2$, from~\eqref{balanced-condition-abelian} we get
$$1-4x=1-4u_2+4s^2>1-4u_2+4u_1^2+4u_2^2=4u_1^2+(1-2u_2)^2\geq 4u_1^2=4y^2,$$
that is, $1-4x>4y^2$ and therefore the Lie algebra
is isomorphic to $\frh_5$.
\end{proof}

It is interesting to point out that the Lie algebras $\frh_2$ and $\frh_4$
have abelian complex structures but it turns out
that none of them admit compatible balanced metric. In contrast, for the Lie
algebra $\frh_5$ we have:

\begin{cor}\label{balanced-h5}
Any abelian complex structure on $\frh_5$ admits balanced Hermitian
metrics.
\end{cor}

\begin{proof}
Abelian complex structures $J$ on $\frh_5$ correspond to $\delta=1$ and
$4y^2<1-4x$. Let us consider $2\,F=i\,(\omega^{1\bar1}+s^2\,\omega^{2\bar2}+t^2\,\omega^{3\bar3})+u\,\omega^{1\bar2}-\bar
u\,\omega^{2\bar1}$ with
$s^2=\frac{1}{2}-x$, $u=y+\frac{i}{2}$ and any non-zero $t^2$. It
is easy to check that $s^2>|u|^2$ and that~\eqref{balanced-condition-abelian} is satisfied.
\end{proof}

The situation is a bit different for the Lie algebra $\frh_3$.
Any complex structure on $\frh_3$ is equivalent to
$\widetilde{J}^+$ or $\widetilde{J}^-$ given by
$$
d\omega^1=d\omega^2=0,\quad d\omega^3=\omega^{1\bar1}\pm\omega^{2\bar2},
$$
but only $\widetilde{J}^-$ admits compatible balanced metrics \cite{U}.

Next we classify, up to equivalence, all the balanced $\widetilde{J}^-$-Hermitian structures on $\frh_3$.

\begin{lemma}\label{reduced_balanced_h3}
Any balanced structure on $(\frh_3,\widetilde{J}^-)$ is equivalent to one and
only one structure in the 1-parameter family
$F=\frac{i}{2}\,(\omega^{1\bar1}+\omega^{2\bar2}+t^2\,\omega^{3\bar3}),\,
t\neq0.$
\end{lemma}

\begin{proof}
First of all, we observe that the balanced condition~\eqref{balanced-condition-abelian} reduces to
$s^2=1$, so the fundamental form of a generic balanced $\widetilde{J}^-$-Hermitian
structure has the following expression:
$$
2F=i\,(\omega^{1\bar1} + \omega^{2\bar2} +
t^2\,\omega^{3\bar3}) +u\,\omega^{1\bar2} - \bar
u\,\omega^{2\bar1},\quad |u|^2<1,\quad t^2>0.
$$

Let us consider the (1,0)-basis $\{\sigma^1,\sigma^2,\sigma^3\}$ given by
$$
\sigma^1 = a_{11}\,\omega^{1} + a_{12}\, \omega^{2},\quad\
\sigma^2 = \overline {a_{12}}\,\omega^{1} + \overline {a_{11}}\, \omega^{2},\quad\
\sigma^3 = a_{33}\,\omega^{3}  ,
$$
where
$$
a_{11}=\left(\frac{1+\sqrt{1-|u|^2}}{2}\right)^{1/2},\quad
a_{12}=\frac{i\bar u}{2}\left(\frac{1+\sqrt{1-|u|^2}}{2}\right)^{-1/2},\quad
a_{33}=\frac{1-|u|^2+\sqrt{1-|u|^2}}{1+\sqrt{1-|u|^2}}.
$$
With respect to this new basis the complex structure equations satisfy
$d\sigma^1=d\sigma^2=0$, $d\sigma^3=\sigma^{1\bar1}-\sigma^{2\bar2}$,
and the fundamental form reduces to
$$2F=i\,(\sigma^{1\bar1} + \sigma^{2\bar2} +
t'^2\,\sigma^{3\bar3})$$
for $t'^2=t^2/|a_{33}|^2$. Now, Lemma~\ref{equivalence-u=0} implies that two structures of this type
corresponding to parameters $t'_1$ and $t'_2$ are equivalent if and only if $t'^2_1=t'^2_2$.
\end{proof}

\subsection{Adapted bases}

One of the main difficulties to study the balanced Hermitian
geometry for nilpotent complex structures $J$ is that the condition (\ref{balanced-condition-nilp}) mixes
the ``metric'' coefficients $s,u$ with the ``complex'' coefficients $b,x,y$ in a
non-trivial way. Next we find an {\it adapted} basis $\{e^1,\ldots,e^6\}$ for any balanced Hermitian structure
in the sense that $\{e^1,\ldots,e^6\}$ is a basis of (real) 1-forms such that the complex structure $J$
and the fundamental 2-form $F$ express in the canonical way
\begin{equation}\label{adapted-basis}
Je^1=-e^2,\ Je^3=-e^4,\ Je^5=-e^6,\quad\quad F=e^{12}+e^{34}+e^{56}.
\end{equation}
It is well-known that such a basis always exists locally, but in the following result
we find an explicit global adapted basis
for any invariant balanced Hermitian structure.

\begin{thrm}\label{general-balanced-real-basis}
Let $(J,F)$ be an invariant balanced Hermitian structure on a $6$-dimensional
nilmanifold~$M$.
Then, there is a basis $\{e^1,\ldots,e^6\}$ of 1-forms on $M$
satisfying \eqref{adapted-basis} and one of the following equations:
\begin{equation}\label{h5real}
de^1 = de^2=de^3=de^4=0,\quad
de^5 = t\,(e^{13} -e^{24}) ,\quad
de^6 = t\,(e^{14} +e^{23}),
\end{equation}
where $t\in\mathbb{R}^*$;
\begin{equation}\label{2-stepreal}
\begin{cases}
\begin{array}{lcl}
de^1 \zzz & = &\zzz de^2=de^3=de^4=0,\\[4pt]
de^5 \zzz & = &\zzz \frac{t}{s}(\rho+b^2)e^{13} -\frac{t}{s}(\rho-b^2) e^{24} ,\\[5pt]
  de^6 \zzz & = &\zzz -2\,t\,(e^{12}-e^{34}) + \frac{t}{s}(\rho-b^2)e^{14}
  + \frac{t}{s}(\rho+b^2)e^{23},
\end{array}
\end{cases}
\end{equation}
where $\rho\in\{0,1\}$, $b\in \mathbb{R}$ and $s,t\in\mathbb{R}^*$;
\begin{equation}\label{realchangeAbelianh5}
\begin{cases}
\begin{array}{lcl}
de^1\zzz & = &\zzz de^2=de^3=de^4=0,\\[6pt]
de^5\zzz & = &\zzz s Y \left[ 2 b^2 u_1 |u|\,(e^{12}-e^{34}) - b^2 t
u_1
|u| Y\,(e^{13} + e^{24}) + 2 \rho s u_1\,(e^{13} - e^{24}) \right.\\[4pt]
&& \quad \left. + 2 s u_2\left((\rho-b^2)e^{14}+(\rho+b^2)e^{23}\right) \right],\\[6pt]
de^6 \zzz & = &\zzz s Y \left[ 2 (2s^2- b^2 u_2)
|u|\,(e^{12}-e^{34})+ b^2 t u_2 |u| Y\,(e^{13}+e^{24})
- 2 \rho s u_2\,(e^{13}-e^{24})\right.\\[4pt]
&& \quad \left. + 2 s
u_1\left((\rho-b^2)e^{14}+(\rho+b^2)e^{23}\right) \right],
\end{array}
\end{cases}
\end{equation}
where $\rho\in\{0,1\}$, $b\in \mathbb{R}$, $t\in\mathbb{R}^*$ and $u\in \mathbb{C}^*$ such that $s^2>|u|^2>0$, and
where $Y=\frac{2 \sqrt{s^2-|u|^2}}{|u|t}$;
\begin{equation}\label{str-eq-Family-I}
\begin{cases}
\begin{array}{lcl}
de^1 \zzz & = &\zzz de^2=de^5=0,\\
de^3 \zzz & = &\zzz \frac{2s}{r}\,e^{15},\\[2pt]  de^4 \zzz & = &\zzz
\frac{2s}{r}\,e^{25},\\[2pt] de^6 \zzz & = &\zzz
\pm\frac{2}{rs}\,(e^{13}+e^{24}),
\end{array}
\end{cases}
\end{equation}
where $r,s\in\mathbb{R}^*$;
\begin{equation}\label{str-eq-Family-II} \begin{cases}
\begin{array}{lcl}
de^1 \zzz & = &\zzz de^2=0,\\[4pt] de^3
\zzz & = &\zzz\frac{s}{rtZ}\,
\left[\pm\frac{t^2}{s^2}\,(e^{13}+e^{24})\pm \frac{t^2}{s^2}(st+Z)\,
(e^{25}-e^{16}) + e^{14}+ \frac{1}{st+Z}\,e^{15}\right],
\\[4pt] de^4 \zzz & = &\zzz\frac{s}{rtZ}\,\left[e^{24}+
\frac{1}{st+Z}\,e^{25}\right],
\\[4pt] de^5\zzz & = &\zzz\frac{-s}{rtZ}\,\left[
(st+Z)\,e^{24}+e^{25}\right],
\\[4pt] de^6 \zzz & = &\zzz\frac{s}{rtZ}\,
\left[\pm \frac{t^2}{s^2}
\frac{1}{st+Z}\,(e^{13}+e^{24})\pm\frac{t^2}{s^2}\,(e^{25}-e^{16}) +
(st+Z)\,e^{14}+ e^{15}\right],
\end{array}
\end{cases}
\end{equation}
where $s,t\in\mathbb{R}^*$ such that $s^2t^2>1$ and where $Z=\sqrt{s^2t^2-1}$.

Furthermore: in case~\eqref{h5real} the complex structure $J$ is complex-parallelizable, i.e. $M$ is the Iwasawa manifold;
for~\eqref{2-stepreal} and~\eqref{realchangeAbelianh5}
the complex structure is nilpotent and the Lie algebra $\frg$ underlying $M$ is isomorphic to $\frh_2$, $\frh_3$, $\frh_4$, $\frh_5$ or $\frh_6$;
abelian complex structures correspond to $\rho=0$ in~\eqref{2-stepreal} and~\eqref{realchangeAbelianh5}, and in such case $\frg\cong\frh_3,\frh_5$
and we can take $b^2=\delta\in\{0,1\}$; finally,
in the remaining cases~\eqref{str-eq-Family-I} and~\eqref{str-eq-Family-II} the Lie algebra underlying $M$ is isomorphic to $\frh_{19}^-$ and the complex structure is non-nilpotent, where
the $\pm$-sign in the equations corresponds to $J=J_0^{\pm}$, respectively.
\end{thrm}

\begin{proof}
Since the balanced Hermitian structure on the nilmanifold $M$ is invariant,
then there is a balanced Hermitian structure $(J,F)$ on the Lie algebra $\frg$ underlying $M$.
We have several possibilities depending on the nilpotency of the complex structure $J$.
If $J$ is complex-parallelizable then by Lemma~\ref{h5_complex_parallelizable}
it suffices to consider the basis $\{e^1,\ldots,e^6\}$ given
by
$$e^1+i\,e^2=\omega^1,\quad e^3+i\,e^4=\omega^2,\quad
e^5+i\,e^6=t\,\omega^3.$$
It is clear that this basis is adapted to $(J,F)$ and the resulting structure equations are~\eqref{h5real}.

Let us suppose now that $J$ is nilpotent but not complex-parallelizable.
We consider two cases in Proposition~\ref{general-balanced}~(ii) depending on
the vanishing of the metric coefficient~$u$.

If $u=0$ then the real basis $\{e^1,\ldots,e^6\}$ given by
$$
e^1+i\,e^2=\omega^1,\quad e^3+i\,e^4=s\,\omega^2,\quad
e^5+i\,e^6=t\,\omega^3,
$$
is a basis adapted to $(J,F)$ and the structure equations become~\eqref{2-stepreal}.

When $u\not=0$, starting from the equations~\eqref{J-nilp} and the balanced condition~\eqref{balanced-condition-nilp},
we consider the $(1,0)$-basis $\{\sigma^1,\sigma^2,\sigma^3\}$ given by
$$
\sigma^1 = u\,\omega^{1} + \frac{i}{2\,\lambda^2}\,
\omega^{2},\quad\quad \sigma^2 = i u\,\omega^{1} -
\frac{1}{2\,\mu^2}\, \omega^2,\quad\quad \sigma^3 = \omega^3, $$
where the coefficients $\lambda^2$ and $\mu^2$ are the roots of the
polynomial $P(X)=s^2 |u|^2 X^2-s^2 X+\frac{1}{4}$, namely
$$\mu^2=\frac{s^2+\sqrt{s^2\,(s^2-|u|^2)}}{2\,s^2|u|^2}>0,\quad
\lambda^2=\frac{1}{|u|^2}-\mu^2>0, \quad \lambda^2\neq \mu^2,\quad
\lambda^2\mu^2=\frac{1}{4s^2|u|^2}.$$
Notice that the roots of $P(X)$ are different, real and strictly
positive because $s^2> |u|^2$. In terms of the new basis the complex structure equations are
\begin{equation*}
\begin{cases}
\begin{array}{lcl}
d\sigma^1 \zzz & = &\zzz d\sigma^2=0,\\[4pt]
d\sigma^3 \zzz & = &\zzz \frac{1}{2s
\sqrt{s^2-|u|^2}}\left[\frac{\rho}{u}\sigma^{12} + (i b^2 \bar{u}
-2s^2)\left(\lambda^2\,\sigma^{1\bar1}
-\mu^2\,\sigma^{2\bar2}\right) + b^2 \bar{u}
\left(\lambda^2\,\sigma^{1\bar2} +
\mu^2\,\sigma^{2\bar1}\right)\right],
\end{array}
\end{cases}
\end{equation*}
and the fundamental form~\eqref{F-nilp} has the simple expression
\begin{equation*}
F=\frac{i}{2}\,(\lambda^2\,\sigma^{1\bar1}+\mu^2\,\sigma^{2\bar2}+t^2\,\sigma^{3\bar3}).
\end{equation*}
Now the real basis $\{e^1,\ldots, e^6\}$ given by
$$
e^1+i\,e^2=\lambda\,\sigma^1,\quad
e^3+i\,e^4=\mu\,\sigma^2,\quad e^5+i\,e^6=t\,\sigma^3,
$$
is clearly a basis adapted to $(J,F)$ and a direct calculation shows that with respect to this basis
the structure equations become~\eqref{realchangeAbelianh5}.

The result for non-nilpotent $J$ follows directly from \cite[Section 3.1]{UV}: starting from Proposition~\ref{general-balanced}~(iii)
it is proved in \cite{UV} that the fundamental form \eqref{F-non-nilp} can be reduced to either $t=1, v=0$ or $v=1$, and these two cases correspond
to equations \eqref{str-eq-Family-I} and \eqref{str-eq-Family-II}, respectively.
\end{proof}

\begin{rmrk}\label{families-equivalence}
{\rm
When the complex structure is non-nilpotent, the classification of balanced Hermitian structures up to equivalence
is obtained in \cite[Theorem 2.10]{UV}, namely: any two balanced Hermitian structures $(J=J_0^{\pm},F_{r,s})$
and $(J'=J_0^{\pm},F_{r',s'})$ given by~\eqref{str-eq-Family-I} are equivalent
if and only if $J=J'$, $r^2=r'^2$ and $s^2=s'^2$; any two balanced Hermitian structures $(J=J_0^{\pm},F_{r,s,t})$
and $(J'=J_0^{\pm},F_{r',s',t'})$ given by~\eqref{str-eq-Family-II} are equivalent
if and only if $J=J'$, $r^2=r'^2$, $s^2=s'^2$ and $t^2=t'^2$;
the structures of family~\eqref{str-eq-Family-I} are not equivalent to the structures of family~\eqref{str-eq-Family-II}

From Lemma~\ref{equivalence-u=0} it follows that, fixed a nilpotent complex structure $J$,
two balanced $J$-Hermitian structures $F_t$ and $F_{t'}$ in the family~\eqref{2-stepreal} are equivalent if and only if $t^2=t'^2$.

A similar result holds for the family~\eqref{h5real} by Lemma~\ref{h5_complex_parallelizable}.
}
\end{rmrk}

\section{Deformation of balanced metrics}\label{deform}

In this section we study on nilmanifolds the weak $\partial\db$-lemma recently introduced in \cite{FY2}
in relation to deformations of balanced metrics.

More precisely, the following definition is given in \cite{FY2}.

\begin{dfn}
{\rm A compact complex $n$-dimensional manifold $M$ satisfies the} $(n-1,n)$-th weak $\partial\db$-lemma
{\rm if for each real form $\varphi$ of type $(n-1,n-1)$ such that $\db\varphi$ is a $\partial$-exact form
there exists a $(n-2,n-1)$-form $\psi$ such that $\db\varphi=i\,\partial\db\psi$.}
\end{dfn}

Fu and Yau prove in \cite{FY2} that given a compact complex $n$-dimensional manifold $M$ with a balanced metric,
if along a small deformation $M_{\lambda}$ of $M$ the $(n-1,n)$-th weak $\partial\db$-lemma is satisfied then
there exists a balanced metric on $M_{\lambda}$ for sufficiently small $\lambda$.

Next we suppose that $M$ is a nilmanifold endowed with an invariant complex structure~$J$ and show when the
weak $\partial\db$-lemma on $(M,J)$ is reduced to the study of the weak $\partial\db$-lemma at the Lie algebra level.

\begin{prop}\label{sym-proc}
Let $M=\nilm$ be a $2n$-dimensional nilmanifold endowed with an invariant complex structure~$J$,
and let $\frg$ be the Lie algebra of $G$.
If $(\frg,J)$ does not satisfy the $(n-1,n)$-th weak $\partial\db$-lemma
then $(M,J)$ does not satisfy the $(n-1,n)$-th weak $\partial\db$-lemma.
\end{prop}

\begin{proof}
The proof is based on the symmetrization process given in \cite{Be} (see also \cite{FG,U}).
Let $\nu=d\tau$ be a volume
element on $M$ induced by a bi-invariant one on the Lie group
$G$ such that, after rescaling, $M$ has volume equal to 1.
Given any covariant $k$-tensor field $T\colon {\mathfrak X}(M)\times
\cdots \times {\frak X}(M) \longrightarrow {\mathcal C}^{\infty}(M)$ on
the nilmanifold $M$, we define a covariant $k$-tensor $T_{\nu}
\colon \frg\times \cdots \times \frg \longrightarrow \mathbb{R}$ on $\frg$ by
$$
T_{\nu}(X_1,\ldots,X_k)=\int_{m\in M} T_m
(X_1\!\mid_m,\ldots,X_k\!\mid_m) \, \nu\ , \quad\quad \mbox{for}\
X_1,\ldots,X_k\in \frg,
$$
where $X_j\!\mid_m$ is the value at the point $m\in M$ of the
projection on $M$ of the left-invariant vector field $X_j$ on the
Lie group $G$. It is clear that
$T_{\nu}=T$ for any tensor field $T$ coming from a
left-invariant one. In~\cite{Be} it is proved that if $T=\alpha$ is a $k$-form on $M$ then
$(d\alpha)_{\nu}=d\alpha_{\nu}$.

Given an invariant complex structure $J$ on $M$ we can extend the symmetrization process to complex forms and it is easy
to see that if $\alpha$ is a form of pure type $(p,q)$ then $\alpha_{\nu}$ is again of pure type $(p,q)$.
Now, for any $(p,q)$-form $\alpha$ on $M$ we have the usual decomposition
$d\alpha=\partial \alpha +\db\alpha$, where $\partial \alpha$ is of type $(p+1,q)$ and $\db\alpha$ of type $(p,q+1)$.
Then, $\partial \alpha_{\nu} +\db \alpha_{\nu}=d\alpha_{\nu}=(d\alpha)_{\nu}=(\partial \alpha)_{\nu} +(\db\alpha)_{\nu}$,
which implies that
$$
(\partial \alpha)_{\nu}=\partial \alpha_{\nu}, \quad (\db \alpha)_{\nu}=\db \alpha_{\nu}.
$$

Suppose that $(\frg,J)$ does not satisfy the $(n-1,n)$-th weak $\partial\db$-lemma, and let
$\varphi$ be a real element in $\bigwedge^{n-1,n-1}(\frg^*)$ such that $\db\varphi=\partial\eta$
for some $\eta\in \bigwedge^{n-2,n}(\frg^*)$
but $\db\varphi \not\in \partial\db \left( \bigwedge^{n-2,n-1}(\frg^*) \right)$.
Therefore, $\varphi$ defines a real $(n-1,n-1)$-form on $M$ such that $\db\varphi=\partial\eta$
but there is no $(n-2,n-1)$-form $\psi$ on $M$ satisfying $\db\varphi=i\,\partial\db\psi$,
because in such case $\psi_{\nu}$ would be an element in $\bigwedge^{n-2,n-1}(\frg^*)$ for which $\db\varphi=i\,\partial\db\, \psi_{\nu}$,
contradicting the fact that $(\frg,J)$ does not satisfy the $(n-1,n)$-th weak $\partial\db$-lemma.
\end{proof}

Let us denote by $H^{p,q}(M,J)$ the Dolbeault cohomology groups
of $(M,J)$ and by $H^{p,q}(\frg,J)$ the cohomology groups of the complex $(\bigwedge^{*,*}(\frg^*),\db)$ at the Lie algebra level.
Conditions under which the natural inclusion
$(\bigwedge^{*,*}(\frg^*),\db) \hookrightarrow (\mathcal{A}^{*,*}(M),\db)$
induces an isomorphism $H^{p,q}(M,J)\cong H^{p,q}(\frg,J)$ in cohomology are investigated in \cite{CF,CFGU2,R}.
In particular, the isomorphism holds for any abelian complex structure $J$.

\begin{rmrk}\label{nomizu}
{\rm
The symmetrization process defines a linear map $\mathcal{A}^{p,q}(M) \longrightarrow \bigwedge^{p,q}(\frg^*)$,
given by $\alpha\mapsto \alpha_\nu$, which commutes with the differentials $\db$. If the natural inclusion
$(\bigwedge^{*,*}(\frg^*),\db) \hookrightarrow (\mathcal{A}^{*,*}(M),\db)$
induces an isomorphism $H^{p,q}(M,J)\cong H^{p,q}(\frg,J)$, then any $\db$-closed $(p,q)$-form $\alpha$ on $M$
is cohomologous to the invariant $(p,q)$-form $\alpha_\nu$ obtained by the symmetrization process.
}
\end{rmrk}

In the next result we find conditions under which the weak $\partial\db$-lemma at the Lie algebra level
implies the weak $\partial\db$-lemma on the nilmanifold.

\begin{prop}\label{non-left-inv}
Let $M=\nilm$ be a $2n$-dimensional nilmanifold endowed with an invariant complex structure~$J$,
and let $\frg$ be the Lie algebra of $G$.
If $(\frg,J)$ satisfies the $(n-1,n)$-th weak $\partial\db$-lemma and $H^{n-2,n}(M,J)\cong H^{n-2,n}(\frg,J)$, then
$(M,J)$ satisfies the $(n-1,n)$-th weak $\partial\db$-lemma.
\end{prop}

\begin{proof}
Let $\varphi$ be a real form of type $(n-1,n-1)$ on $M$ such that $\db\varphi=\partial\eta$ for some
$(n-2,n)$-form $\eta$ on $M$. Since $\db\eta=0$, the form $\eta$ defines a Dolbeault
cohomology class in $H^{n-2,n}(M,J)$. From the isomorphism
$H^{n-2,n}(M,J)\cong H^{n-2,n}(\frg,J)$ and Remark~\ref{nomizu} we get that $\eta=\eta_{\nu}+\db(i\,\psi)$
for some $(n-2,n-1)$-form $\psi$ on $M$. This implies that $\partial\eta=\partial\eta_{\nu}+i\, \partial\db(\psi)$.

Now, from $\db\varphi=\partial\eta$ we get that $\varphi_{\nu}$ is a real element in $\bigwedge^{n-1,n-1}(\frg^*)$ such that $\db\varphi_{\nu}=\partial\eta_{\nu}$. Since $(\frg,J)$ satisfies the $(n-1,n)$-th weak $\partial\db$-lemma, there exists
$\tilde\psi \in \bigwedge^{n-2,n-1}(\frg^*)$ such that $\partial\eta_{\nu}=\db\varphi_{\nu}=i\,\partial\db\tilde\psi$.
Therefore, $\partial\eta=\partial\eta_{\nu}+i\, \partial\db(\psi)=i\, \partial\db(\tilde\psi+\psi)$
and the $(n-1,n)$-th weak $\partial\db$-lemma is satisfied for $(M,J)$.
\end{proof}

Notice that the $(n-1,n)$-th weak $\partial\db$-lemma is satisfied for $(\frg,J)$ if $\partial(\bigwedge^{n-2,n}(\frg^*)) \subset \partial\db(\bigwedge^{n-2,n-1}(\frg^*))$.

\begin{cor}\label{weak-lemma-for-abelian}
Any abelian complex structure satisfies the $(n-1,n)$-th weak $\partial\db$-lemma.
\end{cor}

\begin{proof}
It follows directly from the fact that $\partial(\bigwedge^{n-2,n}(\frg^*))=0$
for any abelian complex structure.
\end{proof}

Next we describe the general behaviour in dimension 6 with respect to the weak $\partial\db$-lemma in the presence of balanced
structures.

\begin{prop}\label{weak-inv}
Let $M$ be a $6$-dimensional nilmanifold endowed with an invariant balanced Hermitian structure $(J,F)$.
Then, the complex manifold $(M,J)$ satisfies the $(2,3)$-th weak $\partial\db$-lemma if and only if $J$ is abelian,
complex-parallelizable or of non-nilpotent type.
\end{prop}

\begin{proof}
The result is known for the Iwasawa manifold, so we suppose next that $J$ is not complex-parallelizable.

Let $\frg$ the Lie algebra underlying $M$ and suppose that $J$ is nilpotent.
By Proposition~\ref{general-balanced} we consider the reduced equations \eqref{J-nilp} and
it is clear that
$$\partial(\bigwedge\!\!\!\!\ ^{1,3}(\frg^*))=\langle \rho\, \omega^{12\bar{1}\bar{2}\bar{3}}\rangle.$$
A direct calculation shows that $\partial\db(\bigwedge^{1,2}(\frg^*))=0$. Now, for the real (2,2)-form $\varphi=\omega^{23\bar{2}\bar{3}}$
we have that $\db\varphi=\omega^{12\bar{1}\bar{2}\bar{3}}$, which implies that the $(2,3)$-th weak $\partial\db$-lemma is not satisfied if
$\rho=1$, i.e. if $J$ is not abelian.
Now, by Proposition~\ref{sym-proc} we get that $(M,J)$ does not satisfy the $(2,3)$-th weak $\partial\db$-lemma if $\rho=1$.

Suppose now that $J$ is non-nilpotent and consider reduced equations as in Proposition~\ref{general-balanced}~(iii).
It is clear that
$\partial(\bigwedge^{1,3}(\frg^*))=\langle \omega^{13\bar{1}\bar{2}\bar{3}}\rangle$
and $\db(\bigwedge^{1,2}(\frg^*))=\langle \omega^{1\bar{1}\bar{2}\bar{3}}, \omega^{2\bar{1}\bar{2}\bar{3}}\rangle$,
which implies that
$$\partial\db(\bigwedge\!\!\!\!\ ^{1,2}(\frg^*))=\langle \omega^{13\bar{1}\bar{2}\bar{3}}\rangle=\partial(\bigwedge\!\!\!\!\ ^{1,3}(\frg^*)).$$
Therefore, the $(2,3)$-th weak $\partial\db$-lemma is satisfied at the Lie algebra level.
Since the Lie algebra $\frg$ is isomorphic to $\frh_{19}^-$ and from \cite{R} the natural inclusion
$(\bigwedge^{*,*}(\frg^*),\db) \hookrightarrow (\mathcal{A}^{*,*}(M),\db)$
induces an isomorphism in cohomology, then the $(2,3)$-th weak $\partial\db$-lemma is satisfied by Proposition~\ref{non-left-inv}.
\end{proof}

The Iwasawa manifold corresponds to the pair $(\frh_5,J_0)$ and it is well-known
that small deformation of the Iwasawa manifold does not admit balanced metric, which implies
that such small deformation does not satisfy the $(2,3)$-th weak $\partial\db$-lemma \cite{FY2}.
In the next example we give, on the nilmanifold $M$
underlying the Iwasawa manifold, an explicit deformation $I_{\lambda}$ of an abelian complex structure $I_0$ having balanced metric such that
the $(2,3)$-th weak $\partial\db$-lemma only holds for $\lambda=0$ but $I_{\lambda}$ admits balanced metric for any~$\lambda$.
This shows that the weak $\partial\db$-lemma is not a necessary condition for the existence of balanced metric along deformation
of the complex structure.

\begin{exm}\label{deform-example}
{\rm
Let us consider $\frh_5$ with basis $e^1,\ldots,e^6$ satisfying $de^1=\cdots=de^4=0$,
$de^5=e^{13}-e^{24}$ and $de^6=e^{14}+e^{23}$. For each $\lambda\in [0,1)$, let us consider the almost complex
structure $I_{\lambda}$ given by
$$
I_{\lambda} e^1=-e^2,\quad I_{\lambda} e^3=-\frac{\lambda+1}{\lambda-1}e^4,\quad
I_{\lambda} e^5=-e^6.
$$
With respect to the basis of (1,0)-forms $\mu^1=e^1+i\, e^2$, $\mu^2=e^3+\frac{\lambda+1}{\lambda-1}i\, e^4$
and $\mu^3=(\lambda+1)(e^5+i\, e^6)$,
the complex structure equations are
$$
d\mu^1=d\mu^2=0,\quad d\mu^3=\lambda\, \mu^{12} +\mu^{1\bar{2}},
$$
which implies the integrability of $I_{\lambda}$. Now it is clear that for $\lambda=0$ the complex structure is abelian and
satisfies the $(2,3)$-th weak $\partial\db$-lemma by Corollary~\ref{weak-lemma-for-abelian}.
Moreover, from Corollary~\ref{balanced-h5} it follows the existence of compatible balanced metric.

When $\lambda\not=0$ we consider the basis of (1,0)-forms $\omega^1=\mu^1$, $\omega^2=\lambda(\mu^2-\mu^1)$ and $\omega^3=\mu^3$,
with respect to which the complex structure equations for $I_{\lambda}$ are
$$
d\omega^1=d\omega^2=0,\quad d\omega^3=\omega^{12} + \omega^{1\bar{1}} + \frac{1}{\lambda} \omega^{1\bar{2}}.
$$
Since these equations are expressed in the form \eqref{J-nilp}, we get by \eqref{balanced-condition-nilp} that $I_{\lambda}$
admits compatible balanced metric if and only if $s^2>|u|^2$ with $s^2 = \frac{\bar u}{\lambda} i$.
Taking $u=\frac{i}{2\lambda}$ and $s^2=\frac{1}{2\lambda^2}$ we have that
the condition $s^2=\frac{1}{2\lambda^2}>\frac{1}{4\lambda^2}=|u|^2$ is satisfied for any $\lambda\in (0,1)$.
In conclusion, for any $\lambda \in (0,1)$ the complex structure $I_{\lambda}$ admits compatible balanced metrics, but
$I_{\lambda}$ does not satisfy the $(2,3)$-th weak $\partial\db$-lemma by Proposition~\ref{weak-inv} because it is nilpotent, but neither complex-parallelizable nor abelian.

Notice that $g_{\lambda}=(e^1)^2+(e^2)^2+\sqrt{\frac{1+\lambda}{1-\lambda}}\, (e^3)^2
+\sqrt{\frac{1+\lambda}{1-\lambda}}\, (e^4)^2+(1+\lambda)(e^5)^2+(1+\lambda)(e^6)^2$
is a balanced $I_{\lambda}$-Hermitian metric for each $\lambda\in [0,1)$.
}
\end{exm}

\section{Holonomy of the Bismut connection}\label{Bismut-hol}

\noindent Bismut proved in \cite{Bis} that any Hermitian structure $(J,F)$ on
a $2n$-dimensional manifold $M$ has a unique Hermitian connection
with torsion $T$ given by $g(X,T(Y,Z))=JdF(X,Y,Z)=-dF(JX,JY,JZ)$,
$g$ being the associated metric. This torsion connection
is known as the {\it Bismut connection} of $(J,F)$ and will be
denoted here by $\nabla$. From now on, we shall identify $T$ with
the 3-form $JdF$. In relation to the Levi-Civita connection
$\nabla^g$ of the Riemannian metric $g$, the Bismut connection is
determined by $\nabla = \nabla^g +\frac{1}{2}T$.

According to \cite{FPS}, the holonomy group of the Bismut connection
associated to any invariant balanced $J$-Hermitian
structure on a nilmanifold $M$ is contained in SU(3).
The aim of this section is to prove that in six dimensions such holonomy group
reduces to a proper subgroup of SU(3) if and only if
the complex structure $J$ is abelian.

In order to prove this result, first we will study explicitly the behaviour of the curvature endomorphisms
of $\nabla$ since they, together with their covariant derivatives, generate the Lie algebra $\mathfrak{hol}(\nabla)$ of the holonomy group
by the well-known Ambrose-Singer theorem.
This approach is also convenient for the applications to the study of the Strominger system in Section~\ref{het-constant-dil}.

The adapted bases found in Theorem~\ref{general-balanced-real-basis} will play a central role.
More concretely, let $\frg$ be a 6-dimensional Lie algebra. Fixed any basis $\{e^1,\ldots,e^6\}$ for the dual $\frg^*$, let us consider the
structure equations
$$
d\, e^k = \sum_{1\leq i<j \leq 6} c_{ij}^k \, e^{i j},\quad\quad
k=1,\ldots,6,
$$
with respect to the basis. Let $g=e^1\otimes e^1 + \cdots +e^6\otimes e^6$ be the inner product
on $\frg$ for which the basis $\{e^k\}^6_{k=1}$ is orthonormal, and
denote by $\{e_1,\ldots,e_6\}$ the dual basis.

Given any linear connection $\nabla$, the connection 1-forms
$\sigma^i_j$ with respect to the basis above are
$$
\sigma^i_j(e_k) = g(\nabla_{e_k}e_j,e_i),
$$
i.e. $\nabla_X e_j = \sigma^1_j(X)\, e_1 +\cdots+ \sigma^6_j(X)\,
e_6$. The curvature 2-forms $\Omega^i_j$ of $\nabla$ are then given in
terms of the connection 1-forms $\sigma^i_j$ by
\begin{equation}\label{curvature}
\Omega^i_j = d \sigma^i_j + \sum_{1\leq k \leq 6}
\sigma^i_k\wedge\sigma^k_j,
\end{equation}
and the curvature endomorphisms $R(e_p,e_q)$ of the connection $\nabla$ are given in terms of the curvature forms~$\Omega^i_j$ by
\begin{equation}\label{endo-curvature}
g(R(e_p,e_q)e_i,e_j) = - \Omega^i_j (e_p,e_q).
\end{equation}

Since $d e^k(e_i,e_j)= -e^k([e_i,e_j])$, the Levi-Civita connection
1-forms $(\sigma^g)^i_j$ of $g$ express in terms of the
structure constants $c_{ij}^k$ by
$$
(\sigma^g)^i_j(e_k) = -\frac12 \left( g(e_i,[e_j,e_k]) - g(e_k,[e_i,e_j]) +
g(e_j,[e_k,e_i]) \right)=\frac12(c^i_{jk}-c^k_{ij}+c^j_{ki}).
$$
Now, let $J$ be a complex structure compatible with $g$ and denote by $F$ the associated fundamental 2-form. Since the Bismut connection $\nabla$ is
given by $\nabla=\nabla^g + \frac12 T$, with torsion $T=JdF$, the Bismut connection 1-forms $\sigma^i_j$ are determined by
\begin{equation}\label{Bismut-1-forms}
\sigma^i_j(e_k)=(\sigma^g)^i_j(e_k) - \frac12 T(e_i,e_j,e_k) =
\frac12(c^i_{jk}-c^k_{ij}+c^j_{ki}) - \frac12 JdF(e_i,e_j,e_k).
\end{equation}

\medskip

Let us suppose next that $\frg$ is the Lie algebra underlying a 6-dimensional nilmanifold $M$ endowed with an invariant Hermitian structure $(J,F)$, and let  $\{e^1,\ldots,e^6\}$ be an adapted basis for $(J,F)$, i.e. satisfying \eqref{adapted-basis}. We can always consider the (3,0)-form $\Psi$ given by
$$
\Psi=(e^1+i\, e^2)\wedge (e^3+i\, e^4)\wedge (e^5+i\, e^6).
$$
Fino, Parton and Salamon prove in \cite{FPS} that the Hermitian structure $(J,F)$ is balanced if and only if $\Psi$ is parallel with respect to the Bismut connection $\nabla$, that is, ${\rm Hol}(\nabla)\subset {\rm SU(3)}$.

We will compute explicitly the Lie algebra $\mathfrak{hol}(\nabla)$ of the holonomy group ${\rm Hol}(\nabla)$ for each invariant balanced Hermitian structure $(J,F)$ on $M$ by using the previous description obtained in Theorem~\ref{general-balanced-real-basis}.
Our main tool is the computation of the curvature endomorphism $R$ and the covariant derivatives, since they generate the Lie algebra $\mathfrak{hol}(\nabla)$. Since we have an adapted basis
$\{e^1,\ldots,e^6\}$ and we know that $\mathfrak{hol}(\nabla)\subset \mathfrak{su}$(3), we will use the following representation
$$
\mathfrak{su}{\rm (3)} \cong \langle \gamma_1,\gamma_2,\gamma_3,\gamma_4,\gamma_5,\gamma_6,\gamma_7,\gamma_8 \rangle
$$
where
\begin{equation}\label{su3-basis}
\begin{array}{lllll}
&\gamma_1=e^{12}-e^{34}, & \gamma_2=e^{13}+e^{24}, & \gamma_3= e^{14}-e^{23}, & \gamma_4=e^{34}-e^{56},\\[7pt]
&\gamma_5=e^{15}+e^{26}, & \gamma_6=e^{16}-e^{25}, & \gamma_7=e^{35}+e^{46}, & \gamma_8=e^{36}-e^{45}.
\end{array}
\end{equation}
Notice that $\gamma_1,\gamma_2,\gamma_3$ generate the Lie subalgebra $\mathfrak{su}(2)$, which will play an important role in the case of abelian complex structures.

Recall that with respect to an adapted bases $\{e^1,\ldots,e^6\}$, the covariant derivative $\nabla_{e_j} \gamma$ of any 2-form $\gamma$
is given by
\begin{equation}\label{cov-deriv}
(\nabla_{e_j} \gamma) (e_p,e_q)= \sum_{k=1}^6 \left( \sigma^k_q(e_j)\, \gamma(e_k,e_p) - \sigma^k_p(e_j)\, \gamma(e_k,e_q) \right),
\quad\ j=1,\ldots, 6.
\end{equation}

\medskip

In order to illustrate the process, we study in the following example the balanced geometry associated to the complex-parallelizable
structure $J_0$, i.e. the standard complex structure on the Iwasawa manifold.

\begin{exm}\label{iwa}
{\rm
The balanced Hermitian geometry associated to the complex-parallelizable structure $J_0$
is described by the structure equations~\eqref{h5real} in Theorem~\ref{general-balanced-real-basis}.
Since the basis $\{e^1,\ldots,e^6\}$ is adapted to the Hermitian structure, by \eqref{h5real} we get that
$dF= t\, e^{136}-t\, e^{145}-t\, e^{235}-t\, e^{246}$ and thus the torsion $T$ is
$$
T=JdF= -t\, e^{135}-t\, e^{146}-t\, e^{236}+t\, e^{245}.
$$
From \eqref{Bismut-1-forms} one has that the non-zero Bismut connection 1-forms $\sigma^i_j$ are the following:
$$
\sigma^1_5=\sigma^2_6= -t\, e^3,\quad\ \sigma^1_6=-\sigma^2_5= -t\, e^4,\quad\ \sigma^3_5=\sigma^4_6= t\, e^1,\quad\
\sigma^3_6=-\sigma^4_5= t\, e^2.
$$
A direct calculation using \eqref{curvature} and \eqref{endo-curvature} gives that the non-zero curvature endomorphisms $R(e_p,e_q)$ of the Bismut connection $\nabla$ are
$$
\begin{array}{lll}
& R(e_1,e_2) = 2t^2 \gamma_4 ,\quad & R(e_1,e_4) = -R(e_2,e_3) = -t^2 \gamma_3 ,\\[4pt]
& R(e_1,e_3) = R(e_2,e_4) = -t^2 \gamma_2 ,\quad & R(e_3,e_4) = 2t^2 \gamma_1 + 2t^2 \gamma_4.
\end{array}
$$
They generate the space $\langle \gamma_1,\gamma_2,\gamma_3,\gamma_4 \rangle$, however, by \eqref{cov-deriv} the covariant derivatives $\nabla_{e_j} \gamma_2$ for $j=1,2,3,4$ are
$$
\begin{array}{lllll}
&\nabla_{e_1} \gamma_2 = -t\, \gamma_5, \quad
& \nabla_{e_2} \gamma_2 = -t\, \gamma_6, \quad
& \nabla_{e_3} \gamma_2 = -t\, \gamma_7, \quad
& \nabla_{e_4} \gamma_2 = -t\, \gamma_8,
\end{array}
$$
and therefore $\mathfrak{hol}(\nabla)=\mathfrak{s}\mathfrak{u}$(3).
}
\end{exm}

For the remaining families of Theorem~\ref{general-balanced-real-basis} the situation is more complicated. We need the following two technical lemmas:

\begin{lemma}\label{nilpotent-u=0}
The curvature endomorphisms $R(e_p,e_q)$ of the Bismut connection $\nabla$ for any structure in the family~\eqref{2-stepreal} are
$$
\begin{array}{llllll}
\frac{s^2}{t^2}\, R(e_1,e_2) \!\!\!&=&\!\!\! -4s^2 \gamma_1 -2 b^2 s\, \gamma_3 +2\rho\,\gamma_4,\quad &\frac{s^2}{t^2}\, R(e_2,e_6) \!\!\!&=&\!\!\! \rho b^2 \gamma_5 - 2 \rho s\, \gamma_7,\\[4pt]
\frac{s^2}{t^2}\, R(e_1,e_3) \!\!\!&=&\!\!\! -(b^4+\rho b^2+\rho) \gamma_2,\quad &\frac{s^2}{t^2}\, R(e_3,e_4) \!\!\!&=&\!\!\! 2(\rho+2s^2) \gamma_1 + 2 b^2 s\, \gamma_3 + 2 \rho \gamma_4,\\[4pt]
\frac{s^2}{t^2}\, R(e_1,e_4) \!\!\!&=&\!\!\! -2(b^2-\rho)s\, \gamma_1 - (b^4-\rho b^2+\rho) \gamma_3,\quad &\frac{s^2}{t^2}\, R(e_3,e_5) \!\!\!&=&\!\!\! \rho b^2 \gamma_7,\\[4pt]
\frac{s^2}{t^2}\, R(e_1,e_5) \!\!\!&=&\!\!\! \rho b^2 \gamma_5,\quad &\frac{s^2}{t^2}\, R(e_3,e_6) \!\!\!&=&\!\!\! 2 \rho s\, \gamma_6 + \rho b^2 \gamma_8,\\[4pt]
\frac{s^2}{t^2}\, R(e_1,e_6) \!\!\!&=&\!\!\! -\rho b^2 \gamma_6 + 2 \rho s\, \gamma_8,\quad &\frac{s^2}{t^2}\, R(e_4,e_5) \!\!\!&=&\!\!\! \rho b^2 \gamma_8,\\[4pt]
\frac{s^2}{t^2}\, R(e_2,e_3) \!\!\!&=&\!\!\! 2(b^2+\rho) s\, \gamma_1 + (b^4+\rho b^2+\rho) \gamma_3,\quad &\frac{s^2}{t^2}\, R(e_4,e_6) \!\!\!&=&\!\!\! - 2 \rho s\, \gamma_5 - \rho b^2 \gamma_7,\\[4pt]
\frac{s^2}{t^2}\, R(e_2,e_4) \!\!\!&=&\!\!\! - (b^4-\rho b^2+\rho) \gamma_2,\quad &\frac{s^2}{t^2}\, R(e_5,e_6) \!\!\!&=&\!\!\! - 2 b^4 \gamma_1 + 4 b^2 s\, \gamma_3.\\[4pt]
\frac{s^2}{t^2}\, R(e_2,e_5) \!\!\!&=&\!\!\! \rho b^2 \gamma_6, &&&
\end{array}
$$
In particular, for $\rho=0$ any $R(e_p,e_q)$ is a linear combination of the following three curvature endomorphisms:
$$
R(e_1,e_2) = -4t^2 \gamma_1 -\frac{2 \delta t^2}{s} \gamma_3, \quad R(e_1,e_3) = -\frac{\delta t^2}{s^2} \gamma_2, \quad
R(e_5,e_6) = - \frac{2 \delta t^2}{s^2} \gamma_1 + \frac{4 \delta t^2}{s} \gamma_3,
$$
where $\delta=b^2\in \{0,1\}$. Moreover, in this case the covariant derivative $\nabla_{e_j} \gamma_i$ is zero for $i=1,2,3$ and $j=1,\ldots,4$, and
$$
\begin{array}{llll}
&\nabla_{e_5} \gamma_1 = \frac{2 \delta t}{s} \gamma_3, \quad
& \nabla_{e_5} \gamma_2 = 0, \quad
& \nabla_{e_5} \gamma_3 = -\frac{2 \delta t}{s} \gamma_1,\\[7pt]
&\nabla_{e_6} \gamma_1 = \frac{2 \delta t}{s} \gamma_2, \quad
& \nabla_{e_6} \gamma_2 = -\frac{2 \delta t}{s} \gamma_1 + 4t \gamma_3, \quad
& \nabla_{e_6} \gamma_3 = -4t \gamma_2.
\end{array}
$$
\end{lemma}

\begin{proof}
It follows from equations \eqref{2-stepreal} that $dF$ is given in terms of the
adapted basis $\{e^1,\ldots,e^6\}$ by
$dF= 2t(e^{125}-e^{345}) + \frac{t}{s} (\rho+b^2) (e^{136}-e^{235})-\frac{t}{s} (\rho-b^2) (e^{145}+e^{246})$, and therefore the torsion is
\begin{equation}\label{torsion-nilpotent-u=0}
T=JdF=-2t(e^{126}-e^{346}) - \frac{t}{s} (\rho-b^2) (e^{135}+e^{236})-\frac{t}{s} (\rho+b^2) (e^{146}-e^{245}).
\end{equation}
By \eqref{Bismut-1-forms} the non-zero Bismut connection 1-forms $\sigma^i_j$ are the following:
$$
\begin{array}{llll}
& \sigma^1_2=-\sigma^3_4 = 2t\, e^6,\quad & \sigma^1_5=\sigma^2_6= -\frac{\rho t}{s}\, e^3,\quad
& \sigma^3_5=\sigma^4_6= \frac{\rho t}{s}\, e^1, \\[6pt]
& \sigma^1_3=\sigma^2_4= -\frac{b^2 t}{s}\, e^5,\quad & \sigma^1_6=-\sigma^2_5= -\frac{\rho t}{s}\, e^4,\quad
& \sigma^3_6=-\sigma^4_5= \frac{\rho t}{s}\, e^2. \\[6pt]
& \sigma^1_4=-\sigma^2_3= \frac{b^2 t}{s}\, e^6, & &
\end{array}
$$
A direct calculation using \eqref{2-stepreal}, \eqref{curvature} and \eqref{endo-curvature} gives the endomorphisms $R(e_p,e_q)$ listed above in terms of the basis \eqref{su3-basis}.
Finally, for $\rho=0$ the covariant derivatives $\nabla_{e_j} \gamma_i$ are easily computed using \eqref{cov-deriv}.
\end{proof}

\begin{lemma}\label{nilpotent-unot=0}
For the balanced Hermitian structures in the family~\eqref{realchangeAbelianh5}, the curvature
endomorphisms $R(e_p,e_q)$ of their Bismut connection $\nabla$ are:
\begin{eqnarray*}
\frac{1}{2|u|^2 s^2 Y^2} R(e_1,e_2) \!\!\!&=&\!\!\! -2(4s^4 +b^4|u|^2-4b^2s^2u_2) \gamma_1
+ b^2t(b^2|u|^2-2s^2u_2)Y \gamma_2 + \frac{4b^2 s^3 u_1}{|u|} \gamma_3 + 4\rho s^2 \gamma_4,\\
\frac{1}{|u|^2 s^2 Y^2} R(e_1,e_3) \!\!\!&=&\!\!\! \frac{2(b^2|u|^2-2s^2u_2)(b^2t|u|Y-2\rho s)}{|u|} \gamma_1
- (b^4t^2|u|^2Y^2-2\rho b^2 s t |u| Y +4\rho s^2) \gamma_2,\\
\frac{1}{4|u|^2 s^4 Y^2} R(e_1,e_4) \!\!\!&=&\!\!\! \frac{2(b^2-\rho)s u_1}{|u|} \gamma_1
- (b^4-\rho b^2 +\rho) \gamma_3,\\
\frac{1}{2s^3Y^2} R(e_1,e_5) \!\!\!&=&\!\!\! -\rho b^2 \left[ (2 s u_2^2 + t u_1^2|u| Y) \gamma_5 + u_1 u_2(2 s- t |u| Y) \gamma_6
+2 u_1 |u| (u_2\gamma_7 + u_1\gamma_8) \right],\\
\frac{1}{2s^3Y^2}  R(e_1,e_6) \!\!\!&=&\!\!\! -\rho \left[ b^2 u_1u_2(2 s - t |u| Y) \gamma_5 + b^2 (2 s u_1^2 + t u_2^2 |u| Y) \gamma_6
-2 |u|(b^2u_2-2s^2) (u_2 \gamma_7 + u_1 \gamma_8) \right],\\
\frac{1}{4|u|^2 s^4 Y^2} R(e_2,e_3) \!\!\!&=&\!\!\! \frac{-2(b^2+\rho)s u_1}{|u|} \gamma_1
+ (b^4+\rho b^2 +\rho) \gamma_3,\\
\frac{1}{|u|^2 s^2 Y^2} R(e_2,e_4) \!\!\!&=&\!\!\! \frac{2(b^2|u|^2-2s^2u_2)(b^2t|u|Y+2\rho s)}{|u|} \gamma_1
- (b^4t^2|u|^2Y^2+2\rho b^2 s t |u| Y +4\rho s^2) \gamma_2,\\
\frac{1}{2s^3Y^2} R(e_2,e_5) \!\!\!&=&\!\!\! \rho b^2 \left[ u_1u_2(2 s - t |u| Y) \gamma_5 -(2 s u_2^2+ t u_1^2|u| Y) \gamma_6
+2 u_1 |u| (u_1\gamma_7 - u_2\gamma_8)\right],\\
\frac{1}{2s^3Y^2} R(e_2,e_6) \!\!\!&=&\!\!\! \rho \left[ b^2(2 s u_1^2 + t u_2^2|u| Y) \gamma_5 - b^2 u_1 u_2(2 s- t |u| Y) \gamma_6
-2 (b^2 u_2-2s^2) |u| (u_1\gamma_7 - u_2\gamma_8) \right] ,\\
\frac{1}{2|u|^2 s^2 Y^2} R(e_3,e_4) \!\!\!&=&\!\!\! 2(2\rho s^2+4s^4 +b^4|u|^2-4b^2s^2u_2) \gamma_1
- b^2t(b^2|u|^2-2s^2u_2)Y \gamma_2 - \frac{4b^2 s^3 u_1}{|u|} \gamma_3 +4\rho s^2 \gamma_4,\\
\frac{1}{2s^3Y^2} R(e_3,e_5) \!\!\!&=&\!\!\! -\rho b^2 \left[ 2u_1 |u| (u_2 \gamma_5 + u_1 \gamma_6)
- (2su_2^2-tu_1^2|u|Y) \gamma_7 - u_1u_2 (2s+t|u|Y) \gamma_8 \right],\\
\frac{1}{2s^3Y^2} R(e_3,e_6) \!\!\!&=&\!\!\! \rho \left[ 2|u|(b^2u_2-2s^2) (u_2 \gamma_5 + u_1 \gamma_6)
+ b^2 u_1u_2(2s+t|u|Y) \gamma_7 + b^2 (2su_1^2-tu_2^2|u|Y) \gamma_8 \right],\\
\frac{1}{2s^3Y^2} R(e_4,e_5) \!\!\!&=&\!\!\! \rho b^2 \left[ 2u_1|u| (u_1 \gamma_5 - u_2 \gamma_6)
- u_1u_2(2s+t|u|Y) \gamma_7 + (2su_2^2-tu_1^2|u|Y) \gamma_8 \right],\\
\frac{1}{2s^3Y^2} R(e_4,e_6) \!\!\!&=&\!\!\! -\rho \left[ 2|u|(b^2u_2-2s^2) (u_1 \gamma_5 - u_2 \gamma_6)
+ b^2 (2su_1^2-tu_2^2|u|Y) \gamma_7 - b^2 u_1u_2 (2s+t|u|Y) \gamma_8 \right],\\
\frac{1}{4|u|^2 s^3 Y^2} R(e_5,e_6) \!\!\!&=&\!\!\! b^4 t|u| Y \gamma_1
+\frac{2b^2(b^2|u|^2-2s^2u_2)}{|u|} \gamma_2 + 2b^2 s t u_1 Y \gamma_3.
\end{eqnarray*}

In particular, for $\rho=0$ any endomorphism $R(e_p,e_q)$ is a linear combination of the following three curvature endomorphisms:
\begin{eqnarray*}
\frac{1}{2|u|^2 s^2 Y^2} R(e_1,e_2) \!\!\!&=&\!\!\! -2(4s^4 +\delta |u|^2-4\delta s^2u_2) \gamma_1
+ \delta t(|u|^2-2s^2u_2)Y \gamma_2 + \frac{4\delta s^3 u_1}{|u|} \gamma_3,\\
\frac{1}{|u|^2 s^2 t Y^3} R(e_1,e_3) \!\!\!&=&\!\!\! 2\delta (|u|^2-2s^2u_2) \gamma_1
- \delta t|u|^2Y \gamma_2,\\
\frac{1}{4|u|^2 s^3 Y^2} R(e_5,e_6) \!\!\!&=&\!\!\! \delta  t|u| Y \gamma_1
+\frac{2\delta(|u|^2-2s^2u_2)}{|u|} \gamma_2 + 2\delta s t u_1 Y \gamma_3,
\end{eqnarray*}
where $\delta=b^2\in\{0,1\}$. Moreover, in this case the covariant derivative $\nabla_{e_j} \gamma_i$
is zero for $i=1,2,3$ and $j=1,\ldots,4$, and
$$
\begin{array}{llll}
& \frac{1}{2sY} \nabla_{e_5} \gamma_1 = 2\delta su_2\, \gamma_2 - \delta tu_1|u|Y\, \gamma_3,\quad
& \frac{1}{2sY} \nabla_{e_6} \gamma_1 = 2\delta su_1\, \gamma_2 + \delta tu_2|u|Y\, \gamma_3,\\[6pt]
& \frac{1}{2sY} \nabla_{e_5} \gamma_2 = -2\delta su_2\, \gamma_1 - 2\delta u_1|u|\, \gamma_3,\quad
& \frac{1}{2sY} \nabla_{e_6} \gamma_2 = -2\delta su_1\, \gamma_1 + 2(\delta u_2-2s^2)|u|\, \gamma_3,\\[6pt]
& \frac{1}{2sY} \nabla_{e_5} \gamma_3 = \delta tu_1|u|Y\, \gamma_1 + 2\delta u_1|u|\, \gamma_2,\quad
& \frac{1}{2sY} \nabla_{e_6} \gamma_3 = -\delta tu_2|u|Y\, \gamma_1 - 2(\delta u_2-2s^2)|u|\, \gamma_2.
\end{array}
$$
\end{lemma}

\begin{proof}
Since the basis $\{e^1,\ldots,e^6\}$ is adapted to the structure, from \eqref{realchangeAbelianh5} we get that the torsion is
\begin{equation}\label{torsion-nilpotent-unot=0}
\begin{array}{lll}
\frac{1}{sY} T= \frac{1}{sY} JdF \!\!\!&=&\!\!\!
2b^2u_1|u|(e^{125}-e^{345}) -2|u|(-2s^2+b^2u_2)(e^{126}-e^{346})\\[6pt]
\!\!\!&&\!\!\! -(2\rho s+b^2t|u|Y) (u_1 e^{135} - u_2 e^{136}) - 2 s(\rho+b^2) (u_2 e^{145} + u_1 e^{146})\\[6pt]
\!\!\!&&\!\!\! -2 s(\rho-b^2) (u_2 e^{235} + u_1 e^{236}) + (2\rho s-b^2t|u|Y) (u_1 e^{245} - u_2 e^{246}).
\end{array}
\end{equation}
Using \eqref{Bismut-1-forms} one has that the non-zero Bismut connection 1-forms $\sigma^i_j$ are the following:
$$
\begin{array}{lll}
& \sigma^1_2=-\sigma^3_4 = -2b^2su_1|u|Y e^5 +2s(-2s^2+b^2u_2)|u|Y e^6,\quad & \sigma^1_6=-\sigma^2_5=
2\rho s^2 Y (u_2 e^3 - u_1 e^4),\\[6pt]
& \sigma^1_3=\sigma^2_4= b^2st|u|Y^2 (u_1 e^5 - u_2 e^6),\quad & \sigma^3_5=\sigma^4_6= 2\rho s^2 Y (u_1 e^1 + u_2 e^2),\\[6pt]
& \sigma^1_4=-\sigma^2_3= 2b^2s^2 Y (u_2 e^5 + u_1 e^6), \quad & \sigma^3_6=-\sigma^4_5= -2\rho s^2 Y (u_2 e^1 - u_1 e^2).\\[6pt]
& \sigma^1_5=\sigma^2_6= -2\rho s^2 Y (u_1 e^3 + u_2 e^4),
\end{array}
$$
A long but direct calculation using \eqref{realchangeAbelianh5}, \eqref{curvature} and \eqref{endo-curvature} gives the endomorphisms $R(e_p,e_q)$ listed above in terms of the basis \eqref{su3-basis}.
Finally, for $\rho=0$ the covariant derivatives $\nabla_{e_j} \gamma_i$ can be computed directly using \eqref{cov-deriv}.
\end{proof}

In the next proposition we describe the Lie algebra $\mathfrak{hol}(\nabla)$ of the holonomy group
of the Bismut connection $\nabla$ when the complex structure $J$ is nilpotent.

\begin{prop}\label{reduc-hol-caso-nilp}
Let $(J,F)$ be an invariant balanced Hermitian structure on a
6-dimensional nilmanifold~$M$ such that $J$ is nilpotent.
Then, the holonomy group of its associated Bismut connection $\nabla$ is equal to {\rm SU(3)}
if and only if the complex structure $J$ is not abelian.

Moreover, if $J$ is abelian then, with respect to the adapted basis satisfying equations~\eqref{2-stepreal} or~\eqref{realchangeAbelianh5} in Theorem~$\ref{general-balanced-real-basis}$ with $\rho=0$, we have:
\begin{itemize}
\item[{\rm (i)}] If the Lie algebra underlying $M$ is isomorphic to
$\frh_5$ then $\mathfrak{hol}(\nabla) \cong \langle \gamma_1,\gamma_2,\gamma_3 \rangle$;
\item[{\rm (ii)}] If the Lie algebra underlying $M$ is isomorphic to
$\frh_3$ then $\mathfrak{hol}(\nabla) \cong \langle \gamma_1 \rangle$.
\end{itemize}
\end{prop}

\begin{proof}
Let us suppose first that $J$ is abelian, i.e. $\rho=0$. From Lemmas~\ref{nilpotent-u=0} and~\ref{nilpotent-unot=0}
it is easy to see that the curvature endomorphisms $R(e_p,e_q)$ and their covariant derivatives of any order lie in the
subspace $\langle \gamma_1,\gamma_2,\gamma_3 \rangle$. We have two possibilities depending on the Lie algebra underlying $M$.
By Proposition~\ref{caracterizacion_h3_h5} the case $\frh_3$ corresponds to $\delta=0$ and it is clear from the previous lemmas that
only $R(e_1,e_2)$ is non-zero and it is proportional to $\gamma_1$, which satisfies $\nabla\gamma_1=0$ and therefore
$\mathfrak{hol}(\nabla) \cong \langle \gamma_1 \rangle$.
On the other hand, the case $\delta=1$ corresponds, by Proposition~\ref{caracterizacion_h3_h5}, to abelian complex structures on $\frh_5$ and it is easy to check from
Lemma~\ref{nilpotent-u=0} that $R(e_1,e_2)$, $R(e_1,e_3)$ and $R(e_5,e_6)$ generate $\langle \gamma_1,\gamma_2,\gamma_3 \rangle$,
so it remains to study the case $\rho=0$ and $\delta=1$ in Lemma~\ref{nilpotent-unot=0}. The determinant of
the matrix whose entries are the components of $R(e_1,e_2)$, $R(e_1,e_3)$ and $R(e_5,e_6)$ in the basis $\{\gamma_1,\gamma_2,\gamma_3\}$ is equal to $\frac{16384\, u_1 (1 + 4 s^2 - 4 u_1^2 - 4 u_2) (s^2 - |u|^2)^{7/2} s^{12}}{|u| t^6}$ and,
since $1 + 4 s^2 - 4 u_1^2 - 4 u_2 > (1-2u_2)^2$ because $s^2-u_2^2>0$, the vanishing of this determinant depends only on the vanishing of $u_1$.
But if $u_1=0$ then a direct calculation shows that $R(e_1,e_3)$ and $R(e_5,e_6)$ generate $\gamma_1$ and $\gamma_2$, and therefore $\gamma_3$ because $[\gamma_1,\gamma_2]=2\gamma_3$. In conclusion, $\mathfrak{hol}(\nabla) \cong \langle \gamma_1,\gamma_2,\gamma_3 \rangle$ when $J$ is an abelian
complex structure on the Lie algebra $\frh_5$.

From now on, let us suppose that the nilpotent complex structure $J$ is not abelian, i.e. $\rho=1$, and we have to prove that
$\mathfrak{hol}(\nabla) \cong \langle \gamma_1,\ldots,\gamma_8 \rangle$. In Example~\ref{iwa} we showed that this holds for the family~\eqref{h5real}. In the case of Lemma~\ref{nilpotent-u=0} it is easy to check
that the curvature endomorphisms $R(e_p,e_q)$ generate the whole space.

For Lemma~\ref{nilpotent-unot=0} we will consider several cases depending on the vanishing of the coefficients $u_1$ and $b$.
Firstly, if $u_1\not=0$ then the determinant of the matrix whose entries are the components of the endomorphisms $\frac{1}{4|u|s^4Y^2} R(e_1,e_4)$ and $\frac{1}{4|u|s^4Y^2} R(e_2,e_3)$ in the basis $\{\gamma_1,\gamma_3\}$ is
$$
\det \begin{pmatrix}
2(b^2-1)s u_1 & -(b^4- b^2 +1)|u|  \\[4pt]
-2(b^2+1)s u_1 & (b^4+ b^2 +1)|u|
\end{pmatrix}
=-4s u_1|u| \not= 0.
$$
From $R(e_1,e_3)$ or $R(e_2,e_4)$ we get $\gamma_2$ because the
components in $\gamma_2$ of these endomorphisms cannot vanish simultaneously, and $\gamma_4$ comes from $R(e_1,e_2)$.

In the case $u_1=0$ we get $\gamma_3$ directly from $R(e_1,e_4)$ or $R(e_2,e_3)$. Again we consider two cases depending on the vanishing of
$b^2 u_2-2s^2$. If $b^2\not= 2s^2/u_2$ then $R(e_1,e_3)$ and $R(e_2,e_4)$ generate $\gamma_1$ and $\gamma_2$, and one also has $\gamma_4$
from $R(e_1,e_2)$. In the case $b^2= 2s^2/u_2$ (which implies $b\not=0$) we get $\gamma_2$ from $R(e_1,e_3)$ or $R(e_2,e_4)$,
$\gamma_1$ from $R(e_5,e_6)$ and then $\gamma_4$ from $R(e_1,e_2)$.

Finally, let us see that the curvature endomorphisms $R(e_p,e_q)$ in Lemma~\ref{nilpotent-unot=0} also generate $\langle \gamma_5,\gamma_6,\gamma_7,\gamma_8 \rangle$.
If $b=0$ then from the curvature endomorphisms $R(e_1,e_6)$, $R(e_2,e_6)$, $R(e_3,e_6)$ and
$R(e_4,e_6)$ it is easy to check that this is true.
For $b\not=0$ we have that the matrix whose entries are the components of $\frac{1}{2b^2s^3Y^2} R(e_1,e_5)$, $\frac{1}{2b^2s^3Y^2} R(e_2,e_5)$, $\frac{1}{2b^2s^3Y^2} R(e_3,e_5)$ and $\frac{1}{2b^2s^3Y^2} R(e_4,e_5)$ in the basis $\{\gamma_5,\gamma_6,\gamma_7,\gamma_8\}$ is
$$
\begin{pmatrix}
-2 s u_2^2 - t u_1^2|u| Y & - u_1 u_2(2 s- t |u| Y) & -2 u_1 u_2 |u| & -2 u_1^2 |u| \\[4pt]
u_1u_2(2 s - t |u| Y) & -2 s u_2^2- t u_1^2|u| Y & 2 u_1^2 |u|  & -2 u_1 u_2 |u| \\[4pt]
-2 u_1 u_2 |u| & -2 u_1^2 |u| & 2su_2^2-tu_1^2|u|Y & u_1u_2 (2s+t|u|Y) \\[4pt]
2 u_1^2 |u| & - 2 u_1 u_2 |u| & -u_1u_2(2s+t|u|Y)  & 2su_2^2-tu_1^2|u|Y
\end{pmatrix}\, ,
$$
whose determinant is equal to $|u|^4 \left( 4 u_1^2 |u|^2 + 4 s^2 u_2^2 +
t^2 u_1^2 |u|^2 Y^2 \right)^2$.
Since it is non-zero, $\gamma_5$, $\gamma_6$, $\gamma_7$ and $\gamma_8$ are generated by these curvature endomorphisms.

In conclusion, if $\rho=1$, i.e. $J$ is nilpotent but non-abelian, then $\mathfrak{hol}(\nabla) \cong \langle \gamma_1,\ldots,\gamma_8 \rangle$
and therefore the holonomy of the Bismut connection is equal to SU(3).
\end{proof}

\begin{rmrk}\label{paralelismo}
{\rm
The reduction of the holonomy of the Bismut connection $\nabla$ to a subgroup of SU(2) in the abelian complex case can also
be derived from the following fact. For $\rho=0$ in the families~\eqref{2-stepreal} and~\eqref{realchangeAbelianh5} we have:
$$
\nabla (e^{12}+e^{34})=0,\quad\  \nabla ((e^1+i\, e^2)\wedge(e^3+i\, e^4))=0,\quad\  \nabla e^5=0,\quad\  \nabla e^6=0.
$$
Moreover, if in addition $\delta=0$ then $\nabla (e^{12})=0$ and $\nabla (e^{34})=0$.
}
\end{rmrk}

It remains to study the holonomy in the case of complex structures of non-nilpotent type.

\begin{prop}\label{reduc-hol-caso-non-nilp}
Let $(J,F)$ be an invariant balanced Hermitian structure on a
6-dimensional nilmanifold~$M$ such that $J$ is non-nilpotent.
Then, the holonomy group of its associated Bismut connection is equal to {\rm SU(3)}.
\end{prop}

\begin{proof}
Since $J$ is non-nilpotent, by Theorem~$\ref{general-balanced-real-basis}$ it is sufficient to study the equations \eqref{str-eq-Family-I} and    \eqref{str-eq-Family-II}. In the first case, we have
that the torsion is
\begin{equation}\label{torsion-fam-I}
T= \pm\frac{2}{rs} e^{136}+\frac{2s}{r} e^{146}
-\frac{2s}{r} e^{236} \pm\frac{2}{rs} e^{246}
\end{equation}
and, by a similar calculation as in the preceding lemmas,
one has
the following curvature endomorphisms of the Bismut connection:
\begin{equation}\label{endo-Bismut-fam-I}
\begin{array}{llll}
& R(e_1,e_2) = -\frac{2s^2}{r^2} \gamma_4 ,\ & R(e_2,e_3) = -R(e_1,e_4) ,\ & R(e_3,e_5) = \frac{s^2}{r^2} \gamma_7,\\[4pt]
& R(e_1,e_3) = \frac{s^4-4}{r^2s^2} \gamma_2 \mp\frac{2}{r^2} \gamma_3 ,\ & R(e_2,e_4) = R(e_1,e_3),\ & R(e_3,e_6) = \pm\frac{2}{r^2} \gamma_7 -\frac{s^2}{r^2} \gamma_8,\\[4pt]
& R(e_1,e_4) = -\frac{s^2}{r^2} \gamma_3,\ & R(e_2,e_5) = \frac{3s^2}{r^2} \gamma_6,\ & R(e_4,e_5) = -\frac{s^2}{r^2} \gamma_8,\\[4pt]
& R(e_1,e_5) = -\frac{3s^2}{r^2} \gamma_5,\ & R(e_2,e_6) = -\frac{s^2}{r^2} \gamma_5 \pm\frac{2}{r^2} \gamma_6,\ & R(e_4,e_6) = -\frac{s^2}{r^2} \gamma_7 \mp\frac{2}{r^2} \gamma_8\\[4pt]
& R(e_1,e_6) = \mp\frac{2}{r^2} \gamma_5 -\frac{s^2}{r^2} \gamma_6,\ & R(e_3,e_4) = -\frac{2s^2}{r^2} (\gamma_1+\gamma_4),\ & R(e_5,e_6) =
R(e_1,e_2) - R(e_3,e_4).
\end{array}
\end{equation}
Thus,
$\mathfrak{hol}(\nabla)\cong \langle \gamma_1,\ldots,\gamma_8 \rangle$.

For the family \eqref{str-eq-Family-II} the torsion $T$ is given by
\begin{eqnarray*}
(rstZ)T \!\!\!&=&\!\!\! s^2(e^{134}-e^{156}) \mp t^2(e^{234}-e^{256})
-(st+Z)(s^2 e^{135} \pm t^2 e^{235})\\
\!\!\!&+&\!\!\! \frac{1}{st+Z}(s^2 e^{146} \pm t^2 e^{246})
+2st(\pm t^2 e^{136}-s^2 e^{236})
\end{eqnarray*}
and one obtains in particular the following curvature endomorphisms of the Bismut connection:
$$
\begin{array}{ll}
& r^2tZ^2 R(e_1,e_4) = \mp 2t \gamma_2 + s^3(\frac{1}{st+Z}-3st) \gamma_3 -s^3 \gamma_5 \pm 2t(st+Z)\gamma_6 ,\\[4pt]
& r^2tZ^2 R(e_1,e_5) = \mp \frac{2t}{st+Z} \gamma_2 -s^3\gamma_3 -s^3(\frac{1}{st+Z}+st) \gamma_5 \pm 2t\gamma_6 ,\\[4pt]
& \frac{r^2Z^2}{s} R(e_2,e_3) = \pm t (st+Z) \gamma_2 + s^3 \gamma_3 \pm t\gamma_6 ,\\[4pt]
& \frac{r^2Z^2}{s} R(e_2,e_6) = \pm t \gamma_2 - s^3 \gamma_5 \pm \frac{t}{st+Z}\gamma_6 .
\end{array}
$$
The determinant of the 4$\times$4 matrix given by the components in $\gamma_2,\gamma_3,\gamma_5$ and $\gamma_6$ of the previous endomorphisms
is equal to $-8\, s^{8}t^{4}$. Since it is non-zero, the endomorphisms $R(e_1,e_4), R(e_1,e_5), R(e_2,e_3)$ and $R(e_2,e_6)$ generate $\gamma_2,\gamma_3,\gamma_5$ and $\gamma_6$.
From the fact that
$$[\gamma_2,\gamma_3]=2\gamma_1, \quad [\gamma_2,\gamma_5]=-\gamma_7,\quad
[\gamma_2,\gamma_6]=-\gamma_8,\quad  [\gamma_5,\gamma_6]=2\gamma_1 +2\gamma_4
$$
we conclude that again $\mathfrak{hol}(\nabla)\cong \langle \gamma_1,\ldots,\gamma_8 \rangle$.

Therefore, if $J$ is non-nilpotent then the holonomy of the Bismut connection always equals SU(3).
\end{proof}

As a consequence of Propositions~\ref{reduc-hol-caso-nilp} and~\ref{reduc-hol-caso-non-nilp} we get:

\begin{thrm}\label{Main1}
Let $(J,F)$ be an invariant balanced Hermitian structure on a
6-dimensional nilmanifold~$M$, and let $\nabla$ be its associated Bismut connection.
Then, ${\rm Hol}(\nabla)={\rm SU(3)}$ if and only if $J$ is not abelian.

Moreover, if the complex structure $J$ is abelian then the holonomy group of the Bismut connection
reduces to a subgroup of {\rm SU(2)},
and it is equal to {\rm SU(2)} if and only if the Lie algebra underlying~$M$ is isomorphic to
the Lie algebra underlying the Iwasawa manifold.
\end{thrm}

\begin{rmrk}\label{jump-hol}
{\rm
Let us consider $\frh_5$ endowed with the balanced Hermitian structures $(I_{\lambda},g_{\lambda})$, $\lambda\in [0,1)$, given in
Example~\ref{deform-example} and let $\nabla_{\lambda}$ denote the associated Bismut connection.
Since $I_\lambda$ is abelian only for $\lambda=0$, by the theorem above
${\rm Hol}(\nabla_0)=$ SU(2) and ${\rm Hol}(\nabla_{\lambda})={\rm SU(3)}$ for any $\lambda\not= 0$.
Notice that $\frh_5$ is the only case where such a `jumping phenomenon' of the Bismut holonomy can occur.
}
\end{rmrk}

The next example shows that Theorem~\ref{Main1} does not hold for abelian complex structures on 6-dimensional compact solvmanifolds.

\begin{exm}\label{abelian-solvable}
{\rm
Let $\frg$ the solvable Lie algebra defined by the equations
$$
de^1=de^2=0,\quad de^3=- e^{13}-e^{24},\quad
de^4=-e^{14}+e^{23},\quad
de^5=e^{15}+e^{26},\quad
de^6=e^{16}-e^{25},
$$
and let us consider $J$ and $F$ defined by \eqref{adapted-basis}. It is easy to check that $J$ is an abelian complex structure
and the pair $(J,F)$ is a balanced Hermitian structure. Since $dF=-2 e^{134}+2 e^{156}$, the torsion is given by
$T=2 e^{234}-2 e^{256}$ and a direct calculation as before shows that $R(e_1,e_2)=0$ and
$$
\begin{array}{lll}
& R(e_1,e_3) = R(e_2,e_4) = -\gamma_2 ,\quad  &  R(e_3,e_4) = 2\, \gamma_1 ,\\[4pt]
& R(e_1,e_4) = -R(e_2,e_3) = -\gamma_3 ,\quad &  R(e_3,e_5) = R(e_4,e_6) = \gamma_7 ,\\[4pt]
& R(e_1,e_5) = R(e_2,e_6) = -\gamma_5 ,\quad  &  R(e_3,e_6) = -R(e_4,e_5) = \gamma_8 ,\\[4pt]
& R(e_1,e_6) = -R(e_2,e_5) = -\gamma_6 ,\quad &  R(e_5,e_6) = 2\,\gamma_1 + 2\,\gamma_4 .
\end{array}
$$
This implies that $\mathfrak{su}(3)\subset \mathfrak{hol}(\nabla)$.
Moreover, the (3,0)-form $\Psi=(e^1+i\, e^2)\wedge (e^3+i\, e^4)\wedge (e^5+i\, e^6)$ is parallel with respect to the Bismut connection, and therefore
$\mathfrak{hol}(\nabla)=\mathfrak{su}(3)$.

The existence of a lattice of maximal rank $\Gamma$ of the simply connected
solvable Lie group $G$ associated to $\frg$ was proved in \cite{Yamada} (see also \cite{deTom}). Therefore, the corresponding
compact solvmanifold has an invariant balanced Hermitian structure $(J,F)$ such that $J$ is abelian and its associated Bismut connection $\nabla$ satisfies ${\rm Hol}(\nabla)={\rm SU(3)}$.
}
\end{exm}

\section{Heterotic supersymmetry with constant
dilaton}\label{het-constant-dil}

\noindent In this section we study the existence of solutions of the Strominger system with respect to the Bismut connection
in the anomaly cancellation condition in the class of abelian complex structures. We show that any invariant balanced metric compatible with
an abelian complex structure provides a solution of the Strominger system.

Since we look for solutions which are invariant, the dilaton function will always be constant.
Recall that a solution of the Strominger system with constant dilaton~\cite{Str} is given by
a compact 6-dimensional manifold $M$ endowed with a Hermitian
SU(3)-structure $(J,F,\Psi)$ satisfying the following system of
equations \cite{Str}:
\begin{enumerate}
\item[(a)] Gravitino equation: the holonomy of the Bismut connection $\nabla$ is contained in
SU(3).
\item[(b)] Dilatino equation with constant dilaton: the dilaton function $\phi$ is constant
and therefore the Lee form $\theta=2d\phi$ is zero, i.e. the Hermitian structure is balanced.
\item[(c)] Gaugino equation: there is a Donaldson-Uhlenbeck-Yau SU(3)-instanton,
i.e. a connection $A$ with curvature 2-forms
$(\Omega^A)^i_j\in \frs\fru(3)$.
\item[(d)] Anomaly cancellation condition:
$dT=2\pi^2 \alpha' \Big(p_1(\nabla)-p_1(A)\Big)$, for $\alpha'>0$.
\end{enumerate}

The instanton $A$ must be non-flat, and $\alpha'$ positive because it is related to the string tension
(for physical interpretation of the solutions of the Strominger system one can see \cite{BBFTY,CCDLMZ,FIUV,FY} and references therein).

In equation (d), $p_1$ denotes the 4-form representing the first Pontrjagin
class of the connection, which is given in terms of the curvature
forms $\Omega^i_j$ of the connection by
$$
p_1= \frac{1}{8\pi^2} {\rm tr}\, \Omega\wedge\Omega
= \frac{1}{8\pi^2} \sum_{1\leq i<j\leq 6} \Omega^i_j\wedge\Omega^i_j.
$$

As we recall in the introduction,
the anomaly cancellation condition could be solved for different metric connections $\nabla$, and we will consider next $\nabla$ as the Bismut connection associated to $(J,F)$.

Let $(J,F)$ be an invariant balanced Hermitian structure on a nilmanifold $M$,
$\{e^1,\ldots,e^6\}$ the adapted basis given in Theorem~\ref{general-balanced-real-basis}
and let us consider the (3,0)-form $\Psi$ defining the
SU(3)-structure given by
$\Psi=(e^1+ie^2)\wedge(e^3+ie^4)\wedge(e^5+ie^6)$.

Notice that in the gaugino equation the curvature 2-forms $(\Omega^A)^i_j\in
{\frak s}{\frak u}(3)$ if and only if
\begin{equation}\label{inst}
(\Omega^A)^i_j(e_1,e_2)+(\Omega^A)^i_j(e_3,e_4)+(\Omega^A)^i_j(e_5,e_6)=0,\quad
(\Omega^A)^i_j(Je_k,Je_l)=(\Omega^A)^i_j(e_k,e_l),\ \ \forall\
i,j,k,l,
\end{equation}
where $\{e_1,\ldots,e_6\}$ is the dual basis of
$\{e^1,\ldots,e^6\}$. We will consider invariant instantons, therefore $A$ satisfies (c) if and
only if each curvature form is a linear combination of the 2-forms $\gamma_1,\ldots,\gamma_8$ given in~\eqref{su3-basis}.

In the next proposition we find SU(3)-instantons for any balanced Hermitian SU(3)-structure
$(J,F,\Psi)$ with abelian $J$.

\begin{prop}\label{instantones-caso-abeliano}
Let $(J,F)$ be any invariant balanced Hermitian structure on a $6$-dimensional nilmanifold~$M$ such that $J$ is abelian.
With respect to the adapted basis $\{e^1,\ldots,e^6\}$ given in Theorem~$\ref{general-balanced-real-basis}$, consider the
{\rm SU(3)}-structure $(J,F,\Psi=(e^1+i\, e^2)\wedge (e^3+i\, e^4)\wedge (e^5+i\, e^6))$.
For each $\lambda\in \mathbb{R}$, let
$A_{\lambda}$ be the {\rm SU(3)}-connection defined by the
connection $1$-forms
\begin{equation}\label{instanton-invariante}
(\sigma^{A_{\lambda}})^1_2=-(\sigma^{A_{\lambda}})^2_1=-(\sigma^{A_{\lambda}})^3_4=(\sigma^{A_{\lambda}})^4_3=\lambda(e^5+e^6),
\end{equation}
and $(\sigma^{A_{\lambda}})^i_j=0$ for $(i,j)\not=(1,2),(2,1),(3,4),(4,3)$. Then,
$A_{\lambda}$ is an {\rm SU(3)}-instanton such that:
\begin{enumerate}
\item[{\rm (i)}] If $(J,F)$ belongs to the family~\eqref{2-stepreal},
then ${\rm tr}\, \Omega^{A_{\lambda}}\wedge \Omega^{A_{\lambda}} =-\frac{8\,t^2}{s^2} (\delta+2s^2) \lambda^2\, e^{1234}$.
\item[{\rm (ii)}] If $(J,F)$ belongs to the family~\eqref{realchangeAbelianh5},
then
$${\rm tr}\, \Omega^{A_{\lambda}}\wedge \Omega^{A_{\lambda}} =-\frac{128\,s^4(s^2-|u|^2)}{t^2}(\delta+2\delta(u_1-u_2)+2s^2) \lambda^2\,e^{1234}.$$
\end{enumerate}
\end{prop}

\begin{proof}
Since $\{e^1,\ldots,e^6\}$ is an adapted basis for the
SU(3)-structure and the connection 1-forms with respect to this
basis satisfy $\sigma^j_i=-\sigma^i_j$ and
$$\sigma^1_3=\sigma^2_4,\quad \sigma^1_4=-\sigma^2_3,\quad
\sigma^1_5=\sigma^2_6,\quad
\sigma^1_6=-\sigma^2_5,\quad\sigma^3_5=\sigma^4_6,\quad
\sigma^3_6=-\sigma^4_5, \quad \sigma^1_2+\sigma^3_4+\sigma^5_6=0,$$
then the connection $A_{\lambda}$ preserves $F$ and $\Psi$, i.e.
it is an SU(3)-connection.

In the case \eqref{2-stepreal} for $\rho=0$, from \eqref{curvature} we get that
$$
\begin{array}{ll}
& (\Omega^{A_{\lambda}})^1_2=-(\Omega^{A_{\lambda}})^2_1=-(\Omega^{A_{\lambda}})^3_4=(\Omega^{A_{\lambda}})^4_3= -2t\lambda\, \gamma_1 +\frac{\delta t\lambda}{s} (\gamma_2 - \gamma_3)
\end{array}
$$
are the only non-zero curvature forms of the connection $A_{\lambda}$.

In the case~\eqref{realchangeAbelianh5} with $\rho=0$, the only non-zero curvature forms are
$$
\begin{array}{ll}
(\Omega^{A_{\lambda}})^1_2=-(\Omega^{A_{\lambda}})^2_1=-(\Omega^{A_{\lambda}})^3_4=(\Omega^{A_{\lambda}})^4_3 & = 2s|u|Y(2s^2+\delta(u_1-u_2))\lambda\, \gamma_1\\[7pt]
& - \delta st|u|Y^2(u_1-u_2)\lambda\, \gamma_2 - 2\delta s^2Y(u_1+u_2)\lambda\, \gamma_3.
\end{array}
$$
Therefore, since the 2-forms $(\Omega^{A_{\lambda}})^i_j$ satisfy
equations~\eqref{inst}, the connection $A_{\lambda}$ is an
SU(3)-instanton in both cases.

Finally, since $\gamma_1\wedge \gamma_1=\gamma_2\wedge \gamma_2=\gamma_3\wedge \gamma_3=-2\, e^{1234}$ and $\gamma_i\wedge \gamma_j=0$
for $1\leq i<j\leq 3$, it is easy to check that the trace of
$\Omega^{A_{\lambda}}\wedge \Omega^{A_{\lambda}}$ is given by (i) or (ii).
\end{proof}

In order to compute the trace of $\Omega\wedge\Omega$ we need to know the curvature forms of the Bismut connection.
However, from \eqref{endo-curvature} and Lemma~\ref{nilpotent-u=0} we get

\begin{lemma}\label{nilpotent-u=0-bis}
The curvature 2-forms $\Omega^i_j$ for the Bismut connection in family~\eqref{2-stepreal} are:
\begin{eqnarray*}
\Omega^1_2 \!\!\!&=&\!\!\! -4t^2(e^{12}-e^{34}) +\frac{2t^2}{s}(\rho-b^2)e^{14} +\frac{2t^2}{s}(\rho+b^2)e^{23}
+\frac{2t^2}{s^2}(\rho\,e^{34} - b^4\,e^{56}),\\
\Omega^1_3 \!\!\!&=&\!\!\! \Omega^2_4 = -\frac{t^2}{s^2}(b^4+\rho\, b^2+\rho) e^{13}
-\frac{t^2}{s^2}(b^4-\rho\, b^2+\rho) e^{24},\\
\Omega^1_4 \!\!\!&=&\!\!\! -\Omega^2_3 = -\frac{2b^2t^2}{s} (e^{12}-e^{34}) -\frac{t^2}{s^2}(b^4-\rho\, b^2+\rho) e^{14}
+\frac{t^2}{s^2}(b^4+\rho\, b^2+\rho) e^{23}
+\frac{4b^2t^2}{s} e^{56},\\
\Omega^1_5 \!\!\!&=&\!\!\! \Omega^2_6 = \frac{\rho\, b^2t^2}{s^2} (e^{15}+e^{26}) - \frac{2 \rho\, t^2}{s} e^{46},\\
\Omega^1_6 \!\!\!&=&\!\!\! -\Omega^2_5 = -\frac{\rho\, b^2t^2}{s^2} (e^{16}-e^{25}) + \frac{2 \rho\,t^2}{s} e^{36},\\
\Omega^3_4 \!\!\!&=&\!\!\! 4 t^2(e^{12}-e^{34}) -\frac{2t^2}{s}(\rho-b^2)e^{14} -\frac{2t^2}{s}(\rho+b^2)e^{23}
+\frac{2t^2}{s^2}(\rho\,e^{12} + b^4\,e^{56}),\\
\Omega^3_5 \!\!\!&=&\!\!\! \Omega^4_6 = - \frac{2 \rho\,t^2}{s} e^{26} + \frac{\rho\,b^2t^2}{s^2} (e^{35}-e^{46}),\\
\Omega^3_6 \!\!\!&=&\!\!\! -\Omega^4_5 = \frac{2 \rho\,t^2}{s} e^{16} + \frac{\rho\,b^2t^2}{s^2} (e^{36}+e^{45}),\\
\Omega^5_6 \!\!\!&=&\!\!\! -\Omega^1_2-\Omega^3_4.
\end{eqnarray*}

Therefore, ${\rm tr}\,\Omega\wedge\Omega = -\frac{8\, t^4}{s^4} (b^8 + \rho b^4 + 4 b^4 s^2 + 2 \rho s^2 + 8 s^4) e^{1234}$.
\end{lemma}

\begin{thrm}\label{solutions-h3}
Let $M$ be a nilmanifold with underlying Lie algebra isomorphic to $\frh_3$.
For any invariant balanced Hermitian structure $(J,F)$ on $M$
there is an invariant {\rm SU(3)}-instanton solving the Strominger system.
Moreover, the Bismut connection associated to $(J,F)$ is of instanton type and therefore
such solutions also solve the heterotic equations of motion.
\end{thrm}

\begin{proof} From Proposition~\ref{caracterizacion_h3_h5} and Lemma~\ref{reduced_balanced_h3}
any invariant balanced Hermitian structure on $M$ is given, up to equivalence,
by equations~\eqref{2-stepreal} with $\rho=0$, $b^2=\delta=0$ and $s^2=1$.
It follows from \eqref{torsion-nilpotent-u=0} that $dT=-8\,t^2 e^{1234}$ and
from Lemma~\ref{nilpotent-u=0-bis} we have ${\rm tr}\,\Omega\wedge\Omega = -64\,t^4 e^{1234}$.
Consider now the instanton $A_{\lambda}$ given in Proposition~\ref{instantones-caso-abeliano}, which by (i) satisfies
${\rm tr}\, \Omega^{A_{\lambda}}\wedge \Omega^{A_{\lambda}} =-16\,t^2 \lambda^2\, e^{1234}$.
The anomaly cancellation condition reduces to solve
$$dT=-8\,t^2 e^{1234}=4\, \alpha'\, t^2(\lambda^2-4\,t^2) e^{1234}=2\pi^2\alpha'\,\left(p_1(\nabla)-p_1(A_{\lambda})\right)$$
for $\alpha'=\frac{2}{4\,t^2-\lambda^2}$ positive. Therefore, it suffices to choose $\lambda$ such that $\lambda^2<4\, t^2$.

Finally, by \cite{I} a solution of the Strominger system is a solution of the heterotic
equations of motion if and only if the connection $\nabla$ in the anomaly cancellation condition is an SU(3)-instanton.
But this is clearly satisfied because by Lemma~\ref{nilpotent-u=0-bis}
the only non-zero curvature forms for the Bismut connection are
$\Omega^1_2 = -\Omega^3_4 = -4t^2\, \gamma_1$.
\end{proof}

To complete the abelian case, we need to consider equations~\eqref{realchangeAbelianh5} with $\rho=0$ and $b^2=\delta=1$.
From Lemma~\ref{nilpotent-unot=0} and relation \eqref{endo-curvature} it follows

\begin{lemma}\label{nilpotent-unot=0-bis}
Let $(J,F)$ be a balanced Hermitian structure in family~\eqref{realchangeAbelianh5} with $\rho=0$ and $b^2=\delta=1$.
The non-zero curvature 2-forms $\Omega^i_j$ of the Bismut connection are:
$$
\begin{array}{ll}
\frac{1}{2s^2|u|^2Y^2}\,\Omega_2^1 = - \frac{1}{2s^2|u|^2Y^2}\,\Omega_4^3 &\!\!\!= -2(4s^4+|u|^2-4s^2u_2)(e^{12}-e^{34}) + tY(|u|^2-2s^2u_2)(e^{13}+e^{24}) \\[10pt]
& +\frac{4s^3u_1}{|u|}(e^{14}-e^{23}) + 2stY|u|\,e^{56},\\[10pt]
\frac{1}{2s^2|u|^2Y^2}\,\Omega_3^1 =  \frac{1}{2s^2|u|^2Y^2}\,\Omega_4^2 =&\!\!\! tY(|u|^2-2s^2u_2)(e^{12}-e^{34}) - \frac{|u|^2t^2Y^2}{2}(e^{13}+e^{24}) + \frac{4s}{|u|}(|u|^2-2s^2u_2)\,e^{56},\\[10pt]
\frac{1}{4s^4|u|Y^2}\,\Omega_4^1 =  -\frac{1}{4s^4|u|Y^2}\,\Omega_3^2 =&\!\!\! 2su_1(e^{12}-e^{34}) - |u|(e^{14}-e^{23}) + 2tu_1|u|Y\,e^{56}.
\end{array}
$$
Therefore, ${\rm tr}\,\Omega\wedge\Omega = -\frac{2048\, s^8 (s^2-|u|^2)^2}{t^4} (1+4s^2 +8s^4 + 4u_1^2 -4u_2 -16s^2u_2+8u_2^2)\, e^{1234}$.
\end{lemma}

Recall that any abelian complex structure on $\frh_5$ admits balanced Hermitian
metrics by Corollary~\ref{balanced-h5}.

\begin{thrm}\label{solutions-h5}
Let $J$ be an abelian complex structure on a nilmanifold $M$ with underlying Lie algebra isomorphic to $\frh_5$.
Then, for any invariant balanced $J$-Hermitian structure on $M$ there exists an invariant {\rm SU(3)}-instanton
solving the Strominger system.
\end{thrm}

\begin{proof} By Proposition~\ref{caracterizacion_h3_h5} the
invariant balanced Hermitian structures on $M$ are given by
equations~\eqref{2-stepreal} or~\eqref{realchangeAbelianh5} with $\rho=0$ and $b^2=\delta=1$.
In the first case, it follows from~\eqref{torsion-nilpotent-u=0} that $dT=-\frac{4\, t^2}{s^2}(1+2s^2) e^{1234}$ and
from Lemma~\ref{nilpotent-u=0-bis} we have ${\rm tr}\,\Omega\wedge\Omega = -\frac{8\, t^4}{s^4}(1+4s^2+8s^4) e^{1234}$.
Consider now the instanton $A_{\lambda}$ given in Proposition~\ref{instantones-caso-abeliano}, which by (i) satisfies
${\rm tr}\, \Omega^{A_{\lambda}}\wedge \Omega^{A_{\lambda}} =-\frac{8\,t^2}{s^2} (1+2s^2) \lambda^2\, e^{1234}$.
Therefore, we need to solve
$$dT=-\frac{4\,t^2}{s^2}(1+2s^2) e^{1234}=\alpha'\,\frac{2\, t^2}{s^4} \left( (1+2s^2)s^2\lambda^2-(1+4s^2+8s^4)t^2 \right) e^{1234}=2\pi^2\alpha'\,\left(p_1(\nabla)-p_1(A_{\lambda})\right)$$
for $\alpha'$ positive. It is sufficient to choose $\lambda$ small enough such that $(1+2s^2)s^2\lambda^2<(1+4s^2+8s^4)t^2$.

Let us consider now the case~\eqref{realchangeAbelianh5} with $\rho=0$ and $b^2=\delta=1$.
From~\eqref{torsion-nilpotent-unot=0} it follows that
$$dT=-\frac{32\, s^2(s^2-|u|^2)}{t^2} \left( s^2+u_1^2+(2s^2-u_2)^2+(s^2-|u|^2) \right) e^{1234}$$
and
from Lemma~\ref{nilpotent-unot=0-bis} we have
$${\rm tr}\,\Omega\wedge\Omega = -\frac{2048\, s^8 (s^2-|u|^2)^2}{t^4} (1+4s^2 +8s^4 + 4u_1^2 -4u_2 -16s^2u_2+8u_2^2)\, e^{1234}.$$
Notice that Lemma~\ref{underlyingLiealgebras}~(ii.3) implies that for any abelian $J$ on $\frh_5$
the condition $4y^2<1-4x$ must be satisfied, where $x=u_2-s^2$ and $y=u_1$.
Therefore, $1+4s^2-4u_2 > 4u_1^2$, which implies that
$$
1+4s^2 +8s^4 + 4u_1^2 -4u_2 -16s^2u_2+8u_2^2 > 8u_1^2 +8(s^2-u_2)^2 \geq 0.
$$
In conclusion, $dT$ and ${\rm tr}\, \Omega\wedge \Omega$ are both a strictly negative multiple of $e^{1234}$.
Consider now the instanton $A_{\lambda}$ given in Proposition~\ref{instantones-caso-abeliano}, which by (ii) satisfies
${\rm tr}\, \Omega^{A_{\lambda}}\wedge \Omega^{A_{\lambda}} =-\frac{128\,s^4(s^2-|u|^2)}{t^2} (1+2s^2+2(u_1-u_2)) \lambda^2\,e^{1234}$.
It is clear that we can choose
$\lambda$ small enough such that $dT=2\pi^2\alpha'\,\left(p_1(\nabla)-p_1(A_{\lambda})\right)$, for $\alpha'$ positive.
\end{proof}

As a consequence of Theorems~\ref{solutions-h3} and~\ref{solutions-h5} any abelian complex structure provides solutions of the Strominger system, more concretely:

\begin{cor}\label{solutions-abelian}
Let $M$ be a nilmanifold endowed with an invariant balanced Hermitian structure $(J,F)$. If $J$ is abelian, then
there is an invariant non-flat {\rm SU(3)}-instanton solving the Strominger system with respect to the Bismut connection in the anomaly
cancellation condition.
Moreover, any such solution solves in addition the heterotic equations of motion if and only if
$M$ is a compact quotient of $H\times \mathbb{R}$, $H$ being the generalized $5$-dimensional Heisenberg group.
\end{cor}

\subsection{More solutions}\label{more-solutions}
As a consequence of the previous study one can also find new solutions of the Strominger system for
complex structures of non-abelian type. For instance, let us consider the family~\eqref{2-stepreal} with $\rho=1$.
From \eqref{torsion-nilpotent-u=0} we have that $dT=-\frac{4\, t^2}{s^2}(1+b^4+2s^2) e^{1234}$, and
by Lemma~\ref{nilpotent-u=0-bis} we get
$${\rm tr}\,\Omega\wedge\Omega = -\frac{8\, t^4}{s^4} (b^8 + b^4 + 4 b^4 s^2 + 2 s^2 + 8 s^4) e^{1234}.$$
Now using the abelian instanton $A$ satisfying ${\rm tr}\,\Omega^A\wedge\Omega^A = -2\, e^{1234}$ given in~\cite{CCDLMZ} we have
$$dT=-\frac{4\, t^2}{s^2}(1+b^4+2s^2) e^{1234}=\frac{\alpha'}{2\, s^4} \left( s^4-4\, t^4 (b^8 + b^4 + 4 b^4 s^2 + 2 s^2 + 8 s^4) \right) e^{1234}=2\pi^2\alpha'\,\left(p_1(\nabla)-p_1(A)\right)$$
with $\alpha'>0$ whenever the metric coefficient $t$ satisfies $4\, t^4 (b^8 + b^4 + 4 b^4 s^2 + 2 s^2 + 8 s^4)>s^4$.

According to Lemma~\ref{underlyingLiealgebras}, for $b^2=1$ the solutions live on a nilmanifold with underlying Lie algebra $\frh_4$ and for $(b^4-1)(b^4-1+4s^2)>0$, resp. $(b^4-1)(b^4-1+4s^2)<0$, the solutions live on a nilmanifold corresponding to the Lie algebra $\frh_5$, resp. $\frh_2$. Also one can prove that the balanced Hermitian structures are not equivalent.

Notice that for $s^2=1$ these solutions were found in \cite[Theorem 6.1]{FIUV}, so
the family above can be
thought as a deformation of such particular solutions.

\medskip

Let $N$ be a nilmanifold with underlying Lie algebra $\frh_{19}^-$, and let us
consider the family of balanced Hermitian
structures $(J_0^{\pm},F)$ given by~\eqref{str-eq-Family-I}.
It follows from~\eqref{torsion-fam-I} that
$$dT=-\frac{8}{r^2}\left(\frac{1}{s^2}\,e^{1234}+s^2\,e^{1256}\right).$$
For each $\tau \in \mathbb{R}$, let
$A_{\tau}$ be the {\rm SU(3)}-connection defined by the
connection $1$-forms
$$(\sigma^{A_{\tau}})^2_3=(\sigma^{A_{\tau}})^2_5=
(\sigma^{A_{\tau}})^4_5=\frac12(\sigma^{A_{\tau}})^5_6=- \tau\,e^6,\qquad
(\sigma^{A_{\tau}})^i_j=\tau\,e^6,$$ for $1\leq i< j\leq 6$ such that $(i,j)\neq (2,3),
(2,5), (4,5), (5,6)$, and $\sigma^j_i=-\sigma^i_j$. By \cite[Proposition~4.1]{UV} the connection
$A_{\tau}$ is an SU(3)-instanton and
\begin{equation}\label{inst-h19}
{\rm tr}\,\Omega^{A_{\tau}}\wedge\Omega^{A_{\tau}}=\frac{-144\,\tau^2}{r^2s^2}\,e^{1234}.
\end{equation}

In the following result we prove that there is a non-flat instanton solving at the same time the
anomaly cancellation conditions for
the Bismut and the Chern connection with respect to the same balanced  Hermitian structure.
To our knowledge, this seems to be the first example with this property.

\begin{prop}\label{simultaneous}
Let $N$ be a nilmanifold with underlying Lie algebra isomorphic to $\frh_{19}^-$. For any invariant complex structure $J$ on $N$
there is a balanced Hermitian structure and a non-flat instanton $A$ solving at the same time the Strominger systems
for the Bismut and the Chern connection.
\end{prop}

\begin{proof}
By \cite{UV} any complex structure $J$ on $\frh_{19}^-$ is equivalent to $J_0^+$ or $J_0^-$, so it suffices to prove the result for $J=J_0^{\pm}$.
We consider the balanced Hermitian structures given in family~\eqref{str-eq-Family-I}. It follows from \cite[Proposition~4.1]{UV}
that for any $r\not=0$ and $s^2\geq 1$, the instanton $A_{\tau}$ with $\tau^2=\frac{s^4-1}{9r^2s^2}$
solves the Strominger system with respect to the Chern connection $\nabla^c$ in the anomaly cancellation condition.

On the other hand, for the Bismut connection $\nabla$, it follows from~\eqref{endo-Bismut-fam-I} that
$$
{\rm tr}\,\Omega\wedge\Omega = \frac{16(s^4-4)}{r^4s^4} e^{1234} - \frac{16\, s^4}{r^4} e^{1256}.
$$
Using~\eqref{inst-h19}, the equation
$
dT=2\pi^2\alpha'\,\left(p_1(\nabla)-p_1(A_{{\widetilde \tau}})\right)
$
has solution if and only if $r\not=0$, $s^2\leq \sqrt{2}$ and ${\widetilde \tau}^2=\frac{2(2-s^4)}{9r^2s^2}$.

Therefore, if $\tau={\widetilde \tau}$ then the corresponding instanton $A$ is non-flat and provides a simultaneous solution.
Notice that $\tau={\widetilde \tau}$ if and only if $s^2=\sqrt{\frac53}$.
\end{proof}

\medskip
\noindent {\bf Acknowledgments.}
Some results of this paper overlap with the second author's Ph. D. thesis
at the University of Zaragoza, under supervision of the first author; we would like to thank M. Fern\'andez,
A. Fino, S. Ivanov, A. Swann and A. Tomassini for their interest and useful comments and suggestions on the work.
This work has been partially supported through Projects MICINN (Spain) MTM2008-06540-C02-02 and MTM2011-28326-C02-01.

\end{document}